\documentclass[12pt]{article}
\usepackage[T1]{fontenc}
\usepackage[latin1]{inputenc}

\usepackage[T1]{fontenc}
\usepackage[latin1]{inputenc}

\usepackage{amsfonts}
\usepackage{amssymb}
\usepackage{amsmath}
\usepackage{eufrak}
\usepackage[dvips]{graphicx}
\usepackage[all]{xy}
\usepackage{latexsym}
\usepackage{verbatim}
\def\Ind#1#2{#1\setbox0=\hbox{$#1x$}\kern\wd0\hbox to 0pt{\hss$#1\mid$\hss}
\lower.9\ht0\hbox to 0pt{\hss$#1\smile$\hss}\kern\wd0}

\def\Notind#1#2{#1\setbox0=\hbox{$#1x$}\kern\wd0\hbox to 0pt{\mathchardef
\nn=12854\hss$#1\nn$\kern1.4\wd0\hss}\hbox to
0pt{\hss$#1\mid$\hss}\lower.9\ht0 \hbox to
0pt{\hss$#1\smile$\hss}\kern\wd0}

\newcommand{\R}{\mathbb R}

\newcommand{\C}{{\EuFrak C}}

\newcommand{\Z}{{\mathbb{Z}}}

\newcommand{\bi}{\begin{itemize}}
\newcommand{\ei}{\end{itemize}}

\newcommand{\tp}{\operatorname{tp}}
\newcommand{\Th}{\operatorname{Th}}

\newcommand{\cl}{\operatorname{cl}}
\newcommand{\id}{\operatorname{id}}
\newcommand{\im}{\operatorname{Im}}
\newcommand{\ext}{\operatorname{ext}}
\newcommand{\defin}{\operatorname{def}}
\newcommand{\extdef}{\operatorname{extdef}}

\DeclareMathOperator{\SL}{SL}
\DeclareMathOperator{\Sp}{Sp}

\newtheorem{theorem}{Theorem}[section]
\newtheorem{lemma}[theorem]{Lemma}
\newtheorem{fact}[theorem]{Fact}

\newtheorem{corollary}[theorem]{Corollary}
\newtheorem{proposition}[theorem]{Proposition}
\newtheorem{definition}[theorem]{Definition}
\newtheorem{remark}[theorem]{Remark}

\newtheorem{question}[theorem]{Question}

\newtheorem{main theorem}{Theorem}
\newtheorem{main question}[main theorem]{Question}
\newtheorem{main conjecture}[main theorem]{Conjecture}
\newtheorem{main corollary}[main theorem]{Corollary}

\title{Generalized Bohr compactification and model-theoretic connected components}
\author{Krzysztof Krupi\'nski\footnote{This research was supported by NCN grant 2012/07/B/ST1/03513.} \hspace{0.5mm} and Anand Pillay\footnote{This research was supported by the EPSRC grant EP/I002294/1, by NSF grant
DMS-1360702, and the MSRI, Berkeley.}}

\date{}

\begin{document}
\maketitle
\begin{abstract}
For a group $G$ first order definable in a structure $M$, we continue the study of the ``definable topological dynamics'' of $G$ (from \cite{GiPePi} for example). The special case when {\em all} subsets of $G$ are definable in the given structure $M$ is simply the usual topological dynamics of the discrete group $G$; in particular, in this case, the words ``externally definable'' and ``definable'' can be removed in the results described below.

Here we consider the mutual interactions of three notions or objects: a certain model-theoretic invariant $G^{*}/(G^{*})^{000}_{M}$ of $G$, which appears to be ``new''  in the classical discrete case and of which we give a direct description in the paper;  the [externally definable] generalized Bohr compactification of $G$; [externally definable] strong amenability.  Among other things, we essentially prove:   (i) The ``new'' invariant  $G^{*}/(G^{*})^{000}_{M}$  lies in between the externally definable generalized Bohr compactification and the definable Bohr compactification, and these all coincide when $G$ is definably strongly amenable and all types in $S_G(M)$ are definable, (ii) the kernel of  the surjective homomorphism from $G^*/(G^*)^{000}_M$ to the
definable Bohr compactification 
%(i.e. the quotient $(G^*)^{00}_M/(G^*)^{000}_M$) 
has naturally the structure of the quotient of a
compact (Hausdorff) group by a dense normal subgroup, and (iii) when $\Th(M)$ is NIP, then $G$ is [externally] definably amenable iff it is externally definably strongly amenable. 

In the situation when all types in $S_G(M)$ are definable, one can just work with the definable (instead of externally definable) objects in the above results.

%In the recent nine years, Newelski has proposed and developed a new, based on topological dynamics, approach to study groups from the model-theoretic perspective. Going further in this direction, for a given group we relate its generalized Bohr compactification with its model-theoretic connected components, and using this, we obtain new information about these components and also about their quotients which are certain model-theoretic invariants of the given group. Our considerations apply in particular to classical discrete topological dynamics by considering the full structure on a given group (i.e. the group structure equipped with predicates for all subsets).  
\end{abstract}

\section{Introduction}

The introduction of methods and ideas from topological dynamics, such as minimal flows and Ellis groups, into the model-theoretic study of definable groups was initiated by Newelski in \cite{Ne1} in the context of the action of a definable group  on its space of types. Among the existing  model-theoretic invariants of a definable group $G$ are  certain  quotients of a saturated elementary extension $G^{*}$ of $G$ by various ``connected components'', which can be also seen as ``quotients'' of the type space, and an important aspect of Newelski's approach \cite{Ne1,Ne2} was to try to relate some of these invariants to dynamical invariants such as the Ellis group.  In the current paper, we try to go deeper in this direction, relating these model-theoretic invariants to the so-called generalized Bohr compactification from \cite{Gl}. As mentioned in the abstract, our model-theoretic set-up subsumes the case when $G$ is simply a discrete group. Our results go both ways. On the one hand, in Theorem \ref{main theorem 1}, we show that the most general model-theoretic invariant $G^{*}/(G^{*})^{000}_{M}$ is a group lying in between the externally definable generalized Bohr compactification and the definable Bohr compactification. In the classical discrete case, $G^{*}/(G^{*})^{000}_{M}$ appears to be a new invariant of $G$, which we will give an intrinsic description of as a universal ``quasihomomorphic'' image of $\beta G$ (see Proposition \ref{prop: description of G/G000}).  On the other hand, the externally definable generalized Bohr compactification gives a new tool for describing $G^{*}/(G^{*})^{000}_{M}$. We should remark that $G^*/(G^*)^{000}_M$ does NOT have the structure of a compact (Hausdorff) group, or even topological group in any meaningful sense, so a lot of previous work has been aimed at trying to understand it as a
mathematical object. And this is one of the main aims and outcomes of the current paper (see Theorem \ref{main theorem 2}). Using \cite{GiKr}, one can also apply our results to get many new examples of discrete groups whose Bohr compactification is trivial whereas the generalized Bohr compactification is non-trivial.

Throughout the paper, we use various notions, theorems and tricks from classical topological dynamics (mostly from \cite{Gl}) adapted to our ``definable'' context. But we would like to emphasize the proofs of our main results are essentially new.

In the rest of this introduction, we will recall some terminology, give background and formulate our main results. More technical definitions will be given in the next section.
%Let us be more precise now (although most of precise definitions will be given in further sections). 

Let $G$ be a group $\emptyset$-definable in a first order structure $M$ (i.e. both the group itself and the graph of the group operation are $\emptyset$-definable sets, that is the sets of realizations in $M$ of some formulas without parameters). By $S_G(M)$ we denote the space of complete types over $M$ containing the formula defining $G$; equivalently, this is the space of ultrafilters 
%containing $G$ 
of definable (with parameters from $M$) subsets of $G$, equipped with the Stone topology. By $S_{G,\ext}(M)$ we denote the space of all externally definable complete types over $M$ containing $G$, that is the space of ultrafilters 
%containing $G$ 
in the Boolean algebra of all externally definable subsets of $G$ (i.e. subsets which are intersections with $G$ of sets definable in arbitrary elementary extensions of $M$). 
%with parameters from arbitrary elementary extensions of $M$).  
The advantage of $S_{G,\ext}(M)$ is that it allows us to develop topological dynamics without any additional assumptions on $M$ (or on the theory of $M$); in order to work smoothly with $S_G(M)$, one would have to assume that all types in $S_G(M)$ are definable (i.e. for every $p \in S_G(M)$ and formula $\varphi(x;y)$ the set of all $m \in M$ such that $\varphi(x,m) \in p$ is definable). Note that if all types in $S_G(M)$ are definable, then $S_G(M)$ coincides with $S_{G,\ext}(M)$. In the case when $M$ is the group $G$ expanded by predicates for all subsets of $G$, all types in $S_G(M)$ are definable, so $S_G(M)$ coincides with $S_{G,\ext}(M)$, and moreover, this is just $\beta G$ -- the Stone-$\check{\mbox{C}}$ech compactification of $G$, that is the space of ultrafilters in the Boolean algebra of all subsets of $G$.

Now, $G$ acts by translations as groups of homeomorphisms of the compact spaces $S_{G,\ext}(M)$, $S_{G}(M)$ and $\beta G$, turning them into $G$-flows (which turn out to be universal in certain categories, see Section \ref{Preliminaries}). Let us focus on $S_{G,\ext}(M)$. The rest of this paragraph is a collection of classical results of Ellis adapted to our context. There is a certain semigroup structure on $S_{G,\ext}(M)$ called the Ellis semigroup; the semigroup operation will be denoted by $*$.  Let ${\cal M}$ be a minimal $G$-subflow of $S_{G,\ext}(M)$. Then ${\cal M}$ is the disjoint union of sets of the form $u*{\cal M}$ with $u$ ranging over all idempotents in ${\cal M}$. Moreover, each $u*{\cal M}$ is a group (with $*$ as the group operation), and the isomorphism type of all these groups is the same and does not depend on the choice of ${\cal M}$. This isomorphism type (or just the group $u * {\cal M}$ for any idempotent $u$) is called the Ellis group of the flow $S_{G,\ext}(M)$.

Let $\C\succ M$ be a monster model extending $M$ (namely, a model which is $\kappa$-saturated (i.e. every type over a set of parameters of cardinality less than $\kappa$ has a realization in this model) and $\kappa$-strongly homogeneous (i.e. every partial elementary map between sets of cardinality less than $\kappa$ extends to an automorphism of the whole model) for a big enough cardinal $\kappa$). By $G^*$ we denote the interpretation of $G$ in $\C$ (i.e. the set of realizations in $\C$ of the formula defining $G$ in $M$). Let $A$ be any small set of parameters from $\C$ (small means of cardinality less than the degree of saturation of $\C$). Recall that an $A$-type-definable set is the set of realizations of a type over $A$ (i.e. it is a (possibly) infinite intersection of $A$-definable sets), and an $A$-invariant set is a set invariant under all automorphisms of the monster model fixing $A$ pointwise.
By a bounded cardinal we will mean a cardinal less than the degree of saturation of $\C$. The following connected components and their quotients play a fundamental role in the study of groups from the model-theoretic perspective: 

\begin{itemize}
\item the intersection of all $A$-definable subgroups of $G^*$ of finite index, denoted by ${G^*}^0_A$,
\item the smallest $A$-type-definable subgroup of $G^*$ of bounded index, denoted by ${G^*}^{00}_A$,
\item the smallest $A$-invariant subgroup of $G^*$ of bounded index, denoted by ${G^*}^{000}_A$.
\end{itemize} 

It is clear that ${G^*}^{000}_A \leq {G^*}^{00}_A \leq {G^*}^{0}_A \leq G^*$, and it is easy to show that all these groups are normal in $G^*$. Sometimes these connected components do not depend on the choice of $A$, e.g. in the so-called NIP theories. In such situations, we skip the parameter set $A$.

The significance of the above components in model theory was discussed in various papers (e.g. see \cite{HrPi,GiNe,Gi}). 
In particular,  they play a fundamental role in the study of stable, simple and NIP groups, and they are precisely related to strong types in different senses \cite{GiNe}. Note also that all quotients such as $G^*/{G^*}^{0}_A$, $G^*/{G^*}^{00}_A$, $G^*/{G^*}^{000}_A$ or ${G^*}^{00}_A/{G^*}^{00}_A$ are certain invariants of the group $G$ (in the sense that they do not depend on the choice of the monster model) and it is desirable to understand their algebraic, topological and possibly some other structure. The components  ${G^*}^0_A$ and ${G^*}^{00}_A$ and the quotients $G^*/{G^*}^{0}_A$ and $G^*/{G^*}^{00}_A$ are well understood. Namely, these quotients are equipped with the logic topology (in which closed sets are the sets with type-definable preimages in $G^*$). Then $G^*/{G^*}^{0}_A$ is a profinite group (in fact, $G^*/{G^*}^{0}_M$ is the definable profinite completion of $G$) and $G^*/{G^*}^{00}_A$ is a compact Hausdorff group (in fact, $G^*/{G^*}^{00}_M$ is the definable Bohr compactification of $G$ which was proved in \cite{GiPePi}). (If one takes as $M$ the group $G$ equipped with predicates for all subsets, then $G^*/{G^*}^{0}_M$ is the profinite completion of $G$ and $G^*/{G^*}^{00}_M$ is the Bohr compactification of $G$). This becomes particularly interesting in $o$-minimal structures due to Pillay's conjecture which describes $G^*/{G^*}^{00}_A$ as a compact Lie group of an appropriate dimension, and so associates with the group $G$ a classical mathematical object ${G^*}/{G^*}^{00}_A$ (see \cite{Pe}). However, ${G^*}^{000}_A$ and the quotients $G^*/{G^*}^{000}_A$ and ${G^*}^{00}_A/{G^*}^{000}_A$ are much more mysterious objects. Until quite recently, it was not known whether always ${G^*}^{00}_A={G^*}^{000}_A$. The first counter-example, which is the universal cover of $\SL_2(\R)$, was found in \cite{CoPi}, and then whole families of counter-examples were found in \cite{GiKr}. Recall that the logic topology on $G^*/{G^*}^{000}_A$ is not necessarily Hausdorff, and on ${G^*}^{00}_A/{G^*}^{000}_A$ it is trivial, so it is rather useless. Thus, one of the 
%%%Krzys1: I replaced "complication" by "complexity".
main goals here is to understand the structure or to learn how to measure the complexity of ${G^*}^{00}_A/{G^*}^{000}_A$. Progress in this direction is done in this paper (see Theorem \ref{main theorem 2} and Corollaries \ref{main corollary 3} and \ref{main corollary 4} below).

Now, we are going to relate the Ellis group $u*{\cal M}$ to connected components and to formulate the main results of this paper. 
 
Newelski \cite{Ne1} found a natural epimorphism, say $\theta$, from $u*{\cal M}$ to $G^*/{G^*}^{00}_M$, and he conjectured that it is an isomorphism. In \cite{GiPePi2}, the authors found a counter-example (namely, the group $G=\SL_2(\R)$ defined in $M=(\R,+,\cdot)$, where $u*{\cal M}\cong\Z_2$ and $G^*/{G^*}^{00}_M$ is trivial, so these groups are not isomorphic), and more examples of this kind were found in \cite{Ja}.
We note in Section \ref{proof of main theorem 1} that there is also a natural epimorphism $f$ from $u*{\cal M}$ to $G^*/{G^*}^{000}_M$ such that $f$ composed with the natural epimorphism $\pi$ from $G^*/{G^*}^{000}_M$ to $G^*/{G^*}^{00}_M$ coincides with $\theta$. This immediately shows that any example where ${G^*}^{000}_M \ne {G^*}^{00}_M$ is a counter-example to Newelski's conjecture (at least in the strong form predicting that $\theta$ is an isomorphism; a weaker form, saying that there is an abstract isomorphism, is also false by the example from \cite{GiPePi2} or by any example from \cite{GiKr} with ${G^*}^{00}_M=G^*\ne {G^*}^{000}_M$). Precise definitions of $f$ and $\theta$ are given in Section \ref{proof of main theorem 1}.

Our goal is to understand better (or to refine) the following sequence of epimorphisms: 

\begin{equation}\label{basic sequence}
\xymatrix{
 u*{\cal M} \ar@{>>}[r]^-{f} & {}{G^*}/{G^*}^{000}_M  \ar@{>>}[r]^-{\pi} & {G^*}/{G^*}^{00}_M}.
\end{equation}

In the context of classical topological dynamics, Glasner (see \cite[Chapter VIII]{Gl}) introduces the generalized Bohr compactification which is meaningful in some situations in which the Bohr compactification is trivial. Analogously, we define an externally definable generalized Bohr compactification of $G$. As in the classical situation, one can introduce a certain topology, called $\tau$-topology, on $u * {\cal M}$ so that if $H(u*{\cal M})$ is the intersection of the $\tau$-closures of $\tau$-neighborhoods of $u$, then $u*{\cal M}/H(u*{\cal M})$ (equipped with the quotient topology induced by $\tau$-topology) becomes a compact Hausdorff group which turns out to be exactly the externally definable generalized Bohr compactification of $G$. However, the proof of this fact (given 
%%%Krzys1: I replaced "slightly" by "somewhat".
in Section \ref{section 2}) in our externally definable context is somewhat different than in the classical case (although it uses various tricks and computations from Glasner's book).

Our first main theorem yields a refinement of the sequence of epimorphisms (\ref{basic sequence}).

\begin{theorem}\label{main theorem 1}
Let $G$ be an arbitrary group $\emptyset$-definable in an arbitrary structure $M$, and let $A\subseteq M$ be any set of parameters.
Let $f: u*{\cal M} \to  G^*/{G^*}^{000}_A$ be the natural epimorphism (defined in Section \ref{proof of main theorem 1}), and suppose that $u*{\cal M}$ is equipped with the $\tau$-topology and $u*{\cal M}/H(u*{\cal M})$ -- with the induced quotient topology. Then:
\begin{enumerate}
\item $f$ is continuous.
\item $H(u*{\cal M}) \leq \ker(f) := f^{-1}(e{G^*}^{000}_A)$.
\item The formula $p*H(u*{\cal M}) \mapsto f(p)$ yields a well-defined continuous epimorphism $\bar f$ from $u*{\cal M}/H(u*{\cal M})$ to $G^*/{G^*}^{000}_A$.
\end{enumerate}
In particular, we get the following sequence of continuous epimorphisms:
\begin{equation}
\xymatrix{
 u*{\cal M} \ar@{>>}[r] &u*{\cal M}/H(u*{\cal M}) \ar@{>>}[r]^-{\bar f}& {}{G^*}/{G^*}^{000}_A  \ar@{>>}[r]^-{\pi} & {G^*}/{G^*}^{00}_A},
\end{equation}
where $u*{\cal M}/H(u*{\cal M})$ is the externally definable generalized Bohr compactification of $G$.
\end{theorem}

Using Theorem \ref{main theorem 1}, we obtain the following property of ${G^*}^{00}_A/{G^*}^{000}_A$.

\begin{theorem}\label{main theorem 2}
For an arbitrary group $G$ which is $\emptyset$-definable in an arbitrary model $M$ and for any set of parameters $A$ the group ${G^*}^{00}_A/{G^*}^{000}_A$ is isomorphic to the quotient of a compact Hausdorff group by a dense subgroup. More precisely, assume $A\subseteq M$, and for $Y:=\ker(\bar f)$ let $\cl_\tau(Y)$ be its closure inside $u*{\cal M}/H(u*{\cal M})$. Then $\bar f$ restricted to $\cl_\tau(Y)$ induces an isomorphism between $\cl_\tau(Y)/Y$ (the quotient of a compact group by a dense subgroup) and ${G^*}^{00}_A/{G^*}^{000}_A$.
\end{theorem}

%%%Krzys: Anand, maybe we should say something about non-commutative geometry in the paragraph below. Can you do that?
The first part of the above theorem has been known to be true in some special situations, e.g. for groups definable in o-minimal expansions of real closed fields \cite{CoPi2} or for some group extensions \cite{GiKr,KrRz}. 
%However, there has been no any known methods so far that could work in general. 
However, there have been no methods so far that could work in general. The goal here is to understand the invariant ${G^*}^{00}_A/{G^*}^{000}_A$ (as was mentioned before, the logic topology on this group is trivial so useless). Theorem \ref{main theorem 2} is interesting in its own right, but it is also possible that it will help us to deduce something new about the Borel cardinality of ${G^*}^{00}_A/{G^*}^{000}_A$ (see \cite{KrPiSo,KaMiSi,KrRz} for definitions) which in a sense measures the complexity of this quotient in descriptive set theoretic terms. 
%Another motivation comes from very recent ideas of Jakub Gismatullin according to which Theorem \ref{main theorem 2} or its variants may lead to new examples 
%%of weakly hyperlinear or even weakly sofic groups, which would be 
%important from the point of view of geometric group theory.
It may also be profitable to view objects such as ${G^*}^{00}_A/{G^*}^{000}_A$ from the point of view of noncommutative geometry, and in fact, Theorem 0.2 gives a representation as a ``typical'' bad quotient of noncommutative geometry. 

We easily get the following corollary of Theorem \ref{main theorem 2}.
\begin{corollary}\label{easy corollary}
Let $A \subseteq M$. Then:\\
i) If the epimorphism $\bar f: u*{\cal M}/H(u*{\cal M}) \to G^*/{G^*}^{000}_A$ is an isomorphism, then ${G^*}^{000}_A={G^*}^{00}_A$.\\
ii) If the epimorphism $f: u*{\cal M} \to G^*/{G^*}^{000}_A$ is an isomorphism, then $H(u*{\cal M})$ is trivial and ${G^*}^{000}_A={G^*}^{00}_A$. 
\end{corollary}

So far we have been working in full generality. In the next result, we need to assume that all types in $S_G(M)$ are definable (still covering the situation when predicates for all subsets of $G$ are added to the language, so still extending the classical discrete topological dynamics). There is a classical notion of a strongly amenable group (see Definition \ref{definition of strong amenability}) and it is claimed in \cite[Chapter IX, Corollary 4.3]{Gl} that for such groups the generalized Bohr compactification is isomorphic to the Bohr compactification. We reprove this theorem in our more general, definable context, showing that the desired isomorphism is exactly the composition of the epimorphism $\bar f:  u*{\cal M}/H(u*{\cal M}) \to G^*/{G^*}^{000}_M$ from Theorem \ref{main theorem 1} with the natural epimorphism $\pi: G^*/{G^*}^{000}_M \to G^*/{G^*}^{00}_M$. Thus, we get the following

\begin{corollary}\label{main corollary 3}
Assume all types in $S_G(M)$ are definable and
suppose $G$ is definably strongly amenable. Then the natural epimorphism $\pi \circ \bar f: u*{\cal M}/H(u*{\cal M}) \to G^*/{G^*}^{00}_M$ is an isomorphism, so ${G^*}^{000}_M={G^*}^{00}_M$. 
\end{corollary}

Now, any nilpotent group $G$  is strongly amenable as a discrete group, so definably strongly amenable for any structure on $G$. Hence:
%Since all nilpotent groups are definably strongly amenable, in particular we have

\begin{corollary}\label{main corollary 4}
Assume all types in $S_G(M)$ are definable and suppose $G$ is nilpotent. Then the natural epimorphism $\pi \circ \bar f : u*{\cal M}/H(u*{\cal M}) \to G^*/{G^*}^{00}_M$ is an isomorphism, so ${G^*}^{000}_M={G^*}^{00}_M$.
\end{corollary}

Since this applies to the situation when predicates for all subsets of $G$ are in the language, Corollary \ref{main corollary 4} yields a partial answer to \cite[Question 4.28]{GiKr} asking for which groups $G$ expanded by predicates for all subsets one has ${G^*}^{000}_G={G^*}^{00}_G$. Recall that by  \cite[Corollary 4.27]{GiKr} (which was obtained by completely different methods), we know that if $G$ is a free group in at least 2 free generators, expanded by predicates for all subsets, then ${G^*}^{000}_G \ne {G^*}^{00}_G$ (and similarly for surface groups of genus at least 2). However, it had not been clear whether this can happen even for abelian groups. Corollary \ref{main corollary 4} shows that that this is impossible in the class of all nilpotent groups. Jakub Gismatullin has told us that he has a different proof of the fact that ${G^*}^{000}_M={G^*}^{00}_M$ in the particular case of abelian groups but in a somewhat more general context (weakening the definability of types assumption).

In all the results above, we used Theorem \ref{main theorem 1} to understand better connected components and their quotients.
There is also an application of Theorem \ref{main theorem 1} which goes the other way around: using some knowledge on connected components together with this theorem, we can produce many examples of discrete groups whose generalized Bohr compactification differs from the Bohr compactification. Namely, working in the situation when predicates for everything are in the language and ${G^*}^{00}_M=G \ne {G^*}^{000}_M$, using  Theorem \ref{main theorem 1} together with the fact that $G^*/{G^*}^{00}_M$ is the Bohr compactification of $G$, we conclude that the Bohr compactification of $G$ is trivial whereas the generalized Bohr compactification is non-trivial. Paper \cite{GiKr} provides many concrete examples where this can be applied, for examples various extensions of symplectic groups $\Sp_{2n}(k)$ over ordered fields $k$.

We also study the relationship between strong amenability and amenability. In topological dynamics, \cite[Chapter III, Theorem 3.1]{Gl} shows that strong amenability implies amenability. Since there are solvable groups which are not strongly amenable, we see that the converse is not true. However, we prove that the counterparts of these notions in a definable context are equivalent for NIP theories.

\begin{theorem}\label{main theorem 5}
If $G$ is a group definable in a NIP theory, then $G$ is externally definably strongly amenable if and only if it is [externally] definably amenable.
\end{theorem} 

We already mentioned that ${G^*}/{G^*}^{00}_M$ is the definable Bohr compactification of $G$.
In the last section, we give a characterization of the externally definable Bohr compactification of $G$ as a quotient of some subgroup of $G$ by a certain component.

\section{Preliminaries}\label{Preliminaries}

In this section, we are going to recall some definitions from topological dynamics, mainly in a modified form to fit the definable or externally definable context. We also make some basic observations. Our main reference for classical topological dynamics is \cite{Gl}, but \cite{Au} may also be helpful.

\subsection{A few basic facts from topological dynamics}

Recall that a $G$-flow is a pair $(G,X)$, where $G$ is a group acting on a compact Hausdorff space by homeomorphisms. We always consider discrete flows, i.e. with no topology on $G$ (or, if one prefers, with the discrete topology on $G$).

%%%Krzys1: I added "the" before "product topology". Is it OK?
\begin{definition}
The Ellis semigroup of the flow $(G,X)$, denoted by $E(X)$, is the closure of the collection of functions $\{\pi_g : g \in G\}$ (where $\pi_g: X \to X$ is given by 
$\pi_g(x)=gx$) in the space $X^X$ equipped with the product topology, with composition as semigroup operation.
\end{definition}

This semigroup operation is continuous in the left coordinate, $E(X)$ is also a $G$-flow, and minimal subflows of $E(X)$ are exactly minimal left ideals with respect to the semigroup structure on $E(X)$. We have the following fundamental fact proved by Ellis.

\begin{fact}\label{Ellis theorem}
Let ${\cal M}$ be a minimal ideal in $E(X)$, and let $J({\cal M})$ be the set of all idempotents in ${\cal M}$. Then the following items hold.\\
i) For any $p \in {\cal M}$, $E(X)p={\cal M}p={\cal M}$.\\
ii) ${\cal M}$ is the disjoint union of sets $u{\cal M}$ with $u$ ranging over $J({\cal M})$.\\
iii) For each $u \in J({\cal M})$, $u{\cal M}$ is a group with the neutral element $u$, where the group operation is the restriction of the semigroup operation on $E(X)$.\\
iv) All the groups $u{\cal M}$ (for $u \in J({\cal M})$) are isomorphic, even when we vary the minimal ideal ${\cal M}$
\end{fact}

The isomorphism type of the groups $u{\cal M}$ (or just any of these groups) from the above fact is called the Ellis group of the flow $X$.

\begin{definition}
A $G$-ambit is a $G$-flow $(G,X,x_0)$ with a distinguished point $x_0 \in X$ such that the orbit $Gx_0$ is dense.
\end{definition} 

With the natural notion of flow homomorphism we have the following fact.

\begin{fact}\label{universality of beta G}
$(G,\beta G, e)$ is the unique up to isomorphism universal $G$-ambit.
\end{fact}

Using this fact, one gets an ``action'' of $\beta G$ on any $G$-flow $(G,X)$, namely for $x \in X$ there is a unique flow homomorphism $h_x: (G,\beta G, e) \to (G,X,x)$, and for $p \in \beta G$ we define $px=h_x(p)$.  More explicitly, this action is given by $px=\lim g_ix$ for any net $(g_i)$ of elements of $G$ converging to $p$ in $\beta G$. In particular, $\beta G$ acts on itself, and denoting this action by $*$, one has $(p*q)x=p(qx)$ for all $p,q \in \beta G$ and $x \in X$. In particular, $*$ is a semigroup operation on $\beta G$ which is continuous on the left and whose restriction to $G \times G$ is the original  group operation on $G$. One easily checks that $(\beta G,*) \cong E(\beta G)$ (by sending $p \in \beta G$ to the function $(x \mapsto px) \in E(\beta G)$). In particular, Fact \ref{Ellis theorem} applies to $(\beta G, *)$ in place of $E(\beta G)$. In a moment, we will recall a very explicit model-theoretic description of $*$.

\subsection{Definable and externally definable context}\label{Definable and externally definable context}

We fix a group $G$ which is $\emptyset$-definable in a first order structure $M$ and a model $N$ which is an $|M|^+$-saturated elementary extension of $M$ (saturation means that every type over less than $|M|^+$ parameters from $N$ has a realization in $N$). We also fix a monster model $\C \succ N$.

By $S_{G,M}(N)$ we denote the space of all types in $S_G(N)$ finitely satisfiable in $M$ (which means that every formula in each such a type has a realization in $M$). Note that $S_{G,M}(N)$ is a closed subset of $S_G(N)$. Then, $S_{G,ext}(M)$ is naturally homeomorphic to $S_{G,M}(N)$, and we identify these spaces. By $S_{G,M}(N)(\C)$ we will denote the set of realizations in $\C$ of all types in $S_{G,M}(N)$ (this is an $N$-type-definable set).
Now, we recall the definition of a definable function from \cite{GiPePi}.

\begin{definition}\label{definable functions}
Let $C$ be a compact space. A map $f: G \to C$ is [externally] definable if for all disjoint, closed subsets $C_1$ and $C_2$ of $C$ the preimages $f^{-1}[C_1]$ and $f^{-1}[C_2]$ can be separated by an [externally] definable subset of $G$.
\end{definition}

For functions defined on the monster model we take another definition, but the lemma below shows that both definitions are compatible.

\begin{definition}\label{Def: ext_def}
Let $C$ be a compact space.\\
i) A function $f: G^* \to C$ is $M$-definable if for every closed $C_1 \subseteq C$ the preimage $f^{-1}[C_1]$ is type-definable over $M$.\\
ii) A function $f: S_{G,M}(N)(\C) \to C$ is externally definable if for every closed $C_1 \subseteq C$ the preimage $f^{-1}[C_1]$ is type-definable over $N$.
\end{definition}

\begin{remark}\label{r,h}
A function $f$ is externally definable according to Definition \ref{Def: ext_def} iff there is a continuous function $h: S_{G,M}(N) \to C$ such that $f=h\circ r$, where $r:S_{G,M}(N)(\C) \to S_{G,M}(N)$ is the obvious map.
\end{remark}

%The argument as in the proof of Lemma 3.2 from \cite{GiPePi} yields

\begin{lemma}\label{prolongation}
1) Definable context:\\
i) If $f: G \to C$ is definable, then it extends uniquely to an $M$-definable function $f^*: G^* \to C$. Moreover, $f^*$ is given by the formula $f^*(a)=\bigcap_{\varphi \in \tp(a/M)} \cl(f[\varphi(M)])$.\\
ii) Conversely, if $f^*:G^* \to C$ is an $M$-definable function, then $f^*|_G: G \to C$ is definable.\\
2) Externally definable context:\\
i) If $f: G \to C$ is externally definable, then it extends uniquely to an externally definable function $f^*: S_{G,M}(N)(\C) \to C$. Moreover, $f^*$ is given by the formula $f^*(a)=\bigcap_{\varphi \in \tp(a/N)} \cl(f[\varphi(M)])$.\\
ii) Conversely, if $f^*:S_{G,M}(N)(\C) \to C$ is an externally definable function, then $f^*|_G: G \to C$ is also externally definable.
\end{lemma}
{\em Proof.} 
(1) is exactly \cite[Lemma 3.2]{GiPePi}, and (2) can be proved in the same fashion. 
\hfill $\square$\\

The definable part of the next definition comes from \cite[Definition 3.5]{GiPePi}.

\begin{definition}\label{definition of definable flows}
i) An [externally] definable $G$-flow is a $G$-flow $(G,X)$ such that for every $x \in X$ the map $f_x: G \to X$ defined by $f_x(g)=gx$ is [externally] definable.\\
ii) An [externally] definable $G$-ambit is an [externally] definable $G$-flow $(G,X,x_0)$ with a distinguished point $x_0$ such that the orbit $Gx_0$ is dense in $X$.
\end{definition}

The proofs of Lemma 3.7(i) and Proposition 3.8 from \cite{GiPePi} can be adapted also to the externally definable context, so we have

\begin{fact}\label{universal definable ambits}
i) If all types in $S_G(M)$ are definable, then $(G,S_G(M),\tp(e/M))$ is the unique up to isomorphism universal definable $G$-ambit.\\ 
ii) $(G,S_{G,M}(N),\tp(e/N))$ is the unique up to isomorphism universal externally definable $G$-ambit.
\end{fact}

As in the case of $\beta G$ in the previous subsection, by universality, we get:
\begin{enumerate}
\item[i)] Assume all types in $S_G(M)$ are definable. Then $S_G(M)$ ``acts'' on every definable $G$-flow $(G,X)$ via $px=\lim g_ix$ for any net $(g_i)$ of elements of $G$ converging to $p$ in $S_G(M)$. In particular, it acts on itself, yielding a semigroup operation $*$ on $S_G(M)$ which is continuous on the left. Then $(S_G(M),*) \cong E(S_G(M))$. In particular, Fact \ref{Ellis theorem} applies to the semigroup $(S_G(M),*)$. The reason why we say that $S_G(M)$ ``acts'' on $(G,X)$ is that $(p *q)x=p(qx)$ for any $p,q \in S_G(M)$ and $x \in X$.
\item[ii)]  $(G,S_{G,M}(N),\tp(e/N))$ ``acts'' on every definable $G$-flow $(G,X)$ via $px=\lim g_ix$ for any net $(g_i)$ of elements of $G$ converging to $p$ in $S_{G,M}(N)$. In particular, it acts on itself, yielding a semigroup operation $*$ on $S_{G,M}(N)$ which is continuous on the left. Then $(S_{G,M}(N),*) \cong E(S_{G,M}(N))$.  In particular, Fact \ref{Ellis theorem} applies to the semigroup $(S_{G,M}(N),*)$. We still have $(p *q)x=p(qx)$ for any $p,q \in S_{G,M}(N)$ and $x \in X$.
\end{enumerate}

Throughout this paper, when we refer to particular results in \cite{Gl},
we often mean (without saying it) the obvious counterparts of the statements
from \cite{Gl} in our [externally] definable context.

The following explicit description of $*$ was established in \cite{Ne2}. This description applies, in particular, to the universal $G$-ambit $(\beta G,*)$ by taking as $M$ the group $G$ equipped with predicates for all subsets of $G$.

\begin{fact}\label{description of *}
i) Assume all types in $S_G(M)$ are definable. Take any $p,q \in S_G(M)$. Then $p*q:=\tp(ab/M)$, where $\tp(b/M)= q$ and $\tp(a/M,b)$ is the unique (by the definability of all types in $S_G(M)$) coheir of $p$ over $M,b$, that is the unique type $\bar p \in S_G(M,b)$ extending $p$ and finitely satisfiable in $M$.\\
ii) Take any $p,q \in S_{G,M}(N)$. Then $p*q:=\tp(ab/N)$, where $\tp(b/N)=q$ and $\tp(a/N,b)$ is the unique extension of $p$ to a complete type over $N,b$ which is finitely satisfiable in $M$.
\end{fact}

The following easy remark is fundamental in many proofs in this paper.

\begin{remark}\label{definability of products}
A product of [externally] definable $G$-flows is an [externally] definable $G$-flow.
\end{remark}
{\em Proof.} 
Let $X_i$, $i \in I$, be [externally] definable $G$-flows. Let $X=\prod_i X_i$, which is naturally a $G$-flow. We will show that it is [externally] definable. Consider any disjoint closed subsets $D_1$ and $D_2$ of $X$ and $x=(x_i)_{i \in I} \in X$. We need to show that the sets $D_1':=\{ g \in G: gx \in D_1\}$ and $D_2':=\{g \in G: gx \in D_2\}$ can be separated by an [externally] definable subset of $G$. Wlog we can assume that $D_1$ and $D_2$ are subbasic closed sets of the form $D_1=A_{1} \times \dots \times A_{n} \times \prod_{i \ne i_1,\dots,i_n} X_i$ and $D_2=B_{1} \times \dots \times B_{n} \times \prod_{i \ne i_1,\dots,i_n} X_i$ for some closed subsets $A_1,B_1 \subseteq X_{i_1}, \dots, A_n,B_n \subseteq X_{i_n}$. 

The assumption that $D_1 \cap D_2=\emptyset$ implies that $A_j \cap B_j=\emptyset$ for some $j \in \{ 1, \dots,n\}$. By the [external] definability of the $G$-flow $X_{i_j}$, we get that $D_1'':=\{ g \in G : gx_{i_j} \in A_j\}$ and $D_2'':=\{ g \in G : gx_{i_j} \in B_j\}$ can be separated by an [externally]  definable subset of $G$. Since clearly $D_1' \subseteq D_1''$ and $D_2' \subseteq D_2''$, the proof is completed. \hfill $\square$\\

The following definition is classical.

\begin{definition}
Two points $x$ and $y$ of a $G$-flow $(G,X)$ are said to be proximal if there exists a net $(g_i)$ in $G$ such that $\lim g_i x = \lim g_i y$. A flow $(G,X)$ is proximal if any two points in it are proximal.
\end{definition}

It is proved in \cite[Chapter II, Proposition 4.2]{Gl} that the unique minimal subflow of the product of all (up to isomorphism) minimal proximal $G$-flows is the unique up to isomorphism universal minimal proximal $G$-flow. Having Remark \ref{definability of products}, Glasner's proof extends to the [externally]  definable context.

\begin{corollary}
There exists a unique up to isomorphism universal minimal [externally] definable proximal $G$-flow, which will be denoted by $(G,\pi_{\defin})$ [respectively, $(G,\pi_{\extdef})$].
\end{corollary}

Recall the definition of amenability and strong amenability in our context.

\begin{definition}\label{definition of strong amenability}
i) $G$ is said to be [externally] definably amenable if there is a left invariant, finitely additive probability measure on [externally] definable subsets of $G$.\\
ii) $G$ is said to be [externally] definably strongly amenable if the only minimal [externally] definable proximal $G$-flow is the trivial one. Equivalently, if the universal  minimal [externally] definable proximal $G$-flow is trivial.
\end{definition}

In classical topological dynamics, it is known that each strongly amenable group is amenable, but the converse is not true. It is also known that all nilpotent groups are strongly amenable (so also [externally] definably strongly amenable).

Recall that a group compactification of a given (discrete) group $G$ is a homomorphism from $G$ to a compact Hausdorff group $K$ with dense image (or just this compact group $K$). (For convenience we will write ``compactification'' instead of ``group compactification''.) There is always a unique up to isomorphism universal compactification of $G$, and it is called the Bohr compactification of $G$. 

By an [externally] definable compactification of $G$ we mean an [externally] definable homomorphism as above.

\begin{remark}
There is a unique universal [externally] definable compactification of $G$, which will be called the [externally] definable Bohr compactification of $G$.
\end{remark}
{\em Proof.} As in the classical case, let $K'$ be the product of all (up to isomorphism) [externally] definable compactifications of $G$, and let $K$ be the closure of the image of the diagonal homomorphism from $G$ (i.e. $K=\cl (\{ (g,g,\dots): g \in G\})$). Clearly $K$ is a compactification of $G$. [External] definability of this compactification follows from Remark \ref{definability of products}. The facts that it is universal and so unique up to isomorphism are obvious. \hfill $\square$\\

In \cite[Proposition 3.4]{GiPePi}, the authors gave the following model-theoretic realization of the definable Bohr compactification.

\begin{fact}\label{Bohr compactification as quotient}
Let $G$ be a group definable in M. Then $G^*/{G^*}^{00}_M$ is the definable Bohr compactification of $G$.
\end{fact}

Almost all the theory recalled and further developed in this paper works in the externally definable case without many differences in comparison with the definable case. However, the last fact is an exception, and the last section of the paper is devoted to a model-theoretic description of the externally definable Bohr compactification of $G$. We will need there the following lemma.

\begin{lemma}\label{partial homomorphism}
Let $f:G \to C$ be an externally definable  homomorphism from $G$ to a compact group $C$. Let $f^*: S_{G,M}(N)(\C) \to C$ be the externally definable map obtained in Lemma \ref{prolongation}, and let $h: S_{G,M}(N) \to C$ be the map provided by Remark \ref{r,h}.\\
i) For any $a_1,\dots,a_n \in S_{G,M}(N)(\C)$ such that $a_1\cdot \ldots \cdot a_n \in S_{G,M}(N)(\C)$ one has $f^*(a_1\cdot \ldots \cdot a_n)=f^*(a_1)\cdot \ldots \cdot f^*(a_n)$.\\
ii) For any $a \in S_{G,M}(N)(\C)$, one has $a^{-1} \in S_{G,M}(N)(\C)$ and $f^*(a^{-1})=f^*(a)^{-1}$.\\
iii) The map $h$ is a semigroup homomorphism.
\end{lemma}
{\em Proof.}
For (i) and (ii) it is enough to apply the argument similar to the one from the proof of \cite[Proposition 3.4]{GiPePi}. In fact, (ii) follows easily from (i).\\[1mm]
iii) Consider any $p,q \in S_{G,M}(N)$. Take $b \models q$ and $a$ such that $\tp(a/N,b)$ is the unique extension of $p$ to a complete type over $N,b$ which is finitely satisfiable in $M$. Then $p*q=\tp(ab/N) \in S_{G,M}(N)$. Therefore, by (ii),  $h(p*q)=f^*(ab)=f^*(a)f^*(b)=h(p)h(q)$. \hfill $\square$\\

One can also think of Fact \ref{Bohr compactification as quotient} as a description of $G^*/{G^*}^{00}_M$ as a universal object in a certain category. Now, we will find such a description for $G^*/{G^*}^{000}_M$.  First, we recall a basic fact on connected components and we fix some notation. 

\begin{fact}\label{G000 as products}
Consider any $A \subseteq M$. Then ${G^*}^{000}_A$ consists of all products of finitely many elements of the form $b^{-1}a$, where $(a,b)$ extends to an infinite $A$-indiscernible sequence, and  ${G^*}^{000}_M$ consists of all products of finitely many elements of the form $b^{-1}a$, where $\tp(a/M)=\tp(b/M)$. In particular, if $\tp(a/M)=\tp(b/M)$, then $b^{-1}a \in {G^*}^{000}_M \subseteq {G^*}^{000}_A$.
\end{fact}

Let $X$ and $Y$ be definable subsets of $G$. Then $X\cdot Y$ will denote
the set of products $x\cdot y$, $x\in X$, $y\in Y$,  viewed first as a definable subset
of $G$, and then as a clopen subset of $S_G(M)$.
For $p, q \in S_G(M)$, we define $p\cdot q$ as the intersection of all clopen sets $X \cdot Y$ for $X\in p$, $Y \in q$,
so a closed subset of $S_G(M)$. Equivalently, $p \cdot q = \{ \tp(a\cdot b/M): \tp(a/M)=p \;\, \mbox{and}\; \tp(b/M)=q\}$.

\begin{proposition}\label{prop: description of G/G000}
Consider surjective maps $f: S_G(M) \to H$, where $H$ is a group, and for
all $p, q\in S_G(M)$, $f(p) \cdot f(q)$ is equal to the common value of  $f(r)$  for all $r\in p\cdot q$.
%and also $ f(p^{-1})=f(p)^{-1})$.
Then there is a (unique) universal such map, and it is precisely the map onto ${G^*}/{G^*}^{000}_M$ given by $\tp(a/M) \mapsto a{G^*}^{000}_M$.

%%%Krzys1: I added "a" before "compact".
When $H$ is a compact Hausdorff group and the surjective map $f:S_G(M) \to H$ is
continuous and $f|_G$ is a homomorphism, then the above condition is automatically
satisfied.
\end{proposition}
{\em Proof.} Uniqueness will be clear once we show that the surjective map $f_0 : S_G(M) \to {G^*}/{G^*}^{000}_M$ given by $f_0(\tp(a/M))= a{G^*}^{000}_M$ is universal. 

First, we check that $f_0$ satisfies the requirement.
The fact that $f_0$ is well-defined (i.e. does not depend on the choice of $a$ realizing a given type) follows from the observation that whenever two elements have the same type over $M$, then they lie in the same coset modulo ${G^*}^{000}_M$ (see Fact \ref{G000 as products}). Now, take  $p,q \in S_G(M)$ and $a \models p$, $b \models q$. Then $f_0(\tp(ab/M))=f_0(\tp(a/M)) \cdot f_0(\tp(b/M))$, which follows immediately from the definition of $f_0$.

%First, we check that $f_0$ satisfies the requirement. So, take $p,q \in S_G(M)$ and $a,a' \models p$ and $b,b' \models q$. We need to check that $f_0(\tp(ab/M))=f_0(\tp(a/M)) \cdot f_0(\tp(b/M))=f_0(\tp(a'b'/M))$. The first equality is obvious from the definition of $f_0$, and the second follows from the basic fact that whenever two elements have the same type over $M$, then they lie in the same coset modulo ${G^*}^{000}_M$ (see Fact \ref{G000 as products}).

Now, consider any $f:S_G(M) \to H$ as in the proposition. Since $f(\tp(e_G/M))=f(\tp(e_G/M) \cdot \tp(e_G/M))=f(\tp(e_G/M))^2$, we see that $f(\tp(e_G/M))=e_H$. So, if for some $a,b \in G^*$ we have $\tp(a/M)=\tp(b/M)$, then $f(\tp(b^{-1}a/M))=f(\tp(b^{-1}/M)) \cdot f(\tp(a/M))=f(\tp(b^{-1}b/M))=e_H$; in particular, $f(p^{-1})=f(p)^{-1}$ for all $p \in S_G(M)$, where $p^{-1}=\tp(a^{-1}/M)$ for any $a \models p$. 

In order to finish the proof, it is enough to show that $h: {G^*}/{G^*}^{000}_M \to H$ given by $h(a{G^*}^{000}_M)=f(\tp(a/M))$ is a well-defined homomorphism. So, take $a,b \in G^*$ such that $a{G^*}^{000}_M=b{G^*}^{000}_M$. By Fact \ref{G000 as products}, this means that $b^{-1}a=\beta_1^{-1}\alpha_1 \cdot \ldots \cdot \beta_n^{-1}\alpha_n$ for some $\alpha_i,\beta_i \in G^*$ such that $\tp(\alpha_i/M)=\tp(\beta_i/M)$ for all i. Therefore, by the last paragraph and the property of $f$, it follows that $f(\tp(a/M))=f(\tp(b/M))$, so we have proved that $h$ is well-defined. The fact that it is a homomorphism follows from the property of $f$.

The additional part of the proposition is an easy exercise using continuity. \hfill $\square$\\

Now, we recall the notion of a group extension from \cite[Chapter VIII]{Gl}.

\begin{definition}
Let $\chi: (G,Y) \to (G,X)$ be a homomorphism of minimal flows. We say that $\chi$ is a group extension if there is a compact Hausdorff group $K$ acting faithfully and continuously on the right on $Y$ in such a way that:\\ 
i) for every $y \in Y$, $\chi^{-1}(\chi(y)) = yK$, and\\
ii) for all $y \in Y$, $g \in G$ and $k \in K$, $(gy)k=g(yk)$.
%A pointed group extension is a homomorpshim $\chi: (G,Y,y_0) \to (G,X,x_0)$ of pointed minimal flows satisfying the above requirements.

In the above situation, we also say that $(G,Y,K)$ is a group extension of $(G,X)$. 
\end{definition}

In particular, each element of $K$ in the above definition is an automorphism of the flow $(G,Y)$. In fact, it is easy to see that $K$ is the group of all automorphisms of $(G,Y)$ preserving the fibers of $\chi$.

We talk about [externally] definable group extension when the flows $(G,X)$ and $(G,Y)$ are [externally] definable.

\begin{remark}
%(i) Let $(G,Y,y_0)$ be a minimal [externally] definable pointed $G$-flow. Then there exists a unique up to isomorphism universal (minimal) [externally] definable pointed group extension $(G,X,K,x_0)$ of $(G,Y,y_0)$.\\
Let $(G,X)$ be a minimal [externally] definable $G$-flow. Then there is a unique up to isomorphism universal (minimal) [externally] definable group extension $(G,Y,K)$ of $(G,X)$. 
\end{remark}
{\em Proof.} We elaborate slightly on the proof from  \cite[page 102]{Gl}, where the author deals with certain pointed $G$-flows. Consider the family $\chi_i: (G,Y_i,K_i) \to (G,X)$, $i \in I$, of all up to isomorphisms [externally] definable group extensions of $(G,X)$.
Fix any point $x_0 \in X$ and choose an idempotent $u \in \beta G$ in some minimal left ideal of $\beta G$ with $ux_0=x_0$ (it exists by \cite[Chapter I, Proposition 3.1(2)]{Gl}). For any $i \in I$ choose $y_i \in  \chi_i^{-1}(x_0)$. We have $\chi_i(uy_i)=u\chi_i(y_i)=ux_0=x_0$, 
%there is $k_i \in K_i$ such that $y_ik_i=uy_i$, so can $y_i'=uy_i$
and so replacing $y_i$ by $uy_i$, we can assume that $uy_i=y_i$ for all $i$.

Let $Z=\prod_i Y_i$ and $z_0 \in Z$ be the sequence $(y_i)$. Finally, put $Y=\cl(Gz_0)$. Since $u$ is an idempotent and clearly $uz_0=z_0$, by \cite[Chapter I, Proposition 3.1]{Gl}, we conclude that $(G,Y)$ is minimal. By Remark \ref{definability of products}, this flow is [externally] definable. Take $K'=\prod_i K_i$ and $K=\{ k \in K': z_0k \in Y\}$. Now, it is easy to check that $(G,Y,K)$ is a universal [externally] definable group extension of $(G,X)$.

The proof of uniqueness on page 102 in \cite{Gl} is not quite complete, because one should work in the category of non-pointed flows.
Consider two universal [externally] definable group extensions $f_1\colon (G,Y_1,K_1) \to (G,X)$ and $f_2\colon (G,Y_2,K_2) \to (G,X)$. Choose any point $y_1 \in Y_1$ and let $x_1=f_1(y_1)$. By the universality of $(G,Y_1,K_1)$ and $(G,Y_2,K_2)$, there are homomorphisms $\varphi \colon (G,Y_1) \to (G,Y_2)$ and $\psi\colon (G,Y_2) \to (G,Y_1)$ satisfying $f_2\circ \varphi=f_1$ and $f_1 \circ \psi =f_2$. Then $\psi(\varphi(y_1)) \in f_1^{-1}(x_1)=y_1K_1$, so, by the minimality of the flow $(G,Y_1)$, we get that the endomorphism $\psi \circ \varphi$ coincides with some $k_1 \in K_1$, so it is an automorphism. Similarly $\varphi \circ \psi$ is an automorphism. Thus, $\varphi$ is an isomorphism. \hfill $\square$\\

Now, we recall one of the central notions for this paper (see \cite[Chapter VIII]{Gl}).

\begin{definition}
The universal [externally] definable group extension $(G,\pi_{\defin}^\#,\Sigma_{\defin}(G))$ [respectively, $(G,\pi_{\extdef}^\#,\Sigma_{\extdef}(G))$] of the universal minimal [externally] definable proximal $G$-flow $(G,\pi_{def})$ [respectively, $(G,\pi_{\extdef})$] will be called the universal [externally] definable compactification flow of $G$, and $\Sigma_{\defin}(G)$ [respectively,   $\Sigma_{\extdef}(G)$]-- the generalized [externally]  definable Bohr compactification of $G$.
\end{definition}

Recall that a minimal $G$-flow $(G,X)$ is said to be regular if  for any $x,y \in X$ there is an automorphism $f$ of $X$ such that $f(x)$ and $y$ are proximal. A regular $G$-flow $(G,X)$ is called a compactification flow of $G$ if the group of all automorphisms of $(G,X)$ is a compact Hausdorff topological group (with the pointwise convergence topology), and this compact group of automorphisms is called a generalized compactification of $G$ (note that it is not necessarily a compactification of $G$). It is proved in \cite[Chapter VIII, Proposition 3.1]{Gl} that $(G,X)$ is a compactification flow of $G$ if and only if it is a group extension of a proximal flow, and it follows from the proof that the group $K$ associated with this group extension is exactly the group of all automorphisms of $(G,X)$. In particular, $\Sigma_{\defin}(G)$ [respectively,  $\Sigma_{\extdef}(G)$] is the compact group of all automorphisms of  $(G,\pi_{\defin}^\#)$ [respectively, of $(G,\pi_{\extdef}^\#)$].

In this paper, we will usually work in $(S_{G,\ext}(M),*)$. So ${\cal M}$ will be a minimal left ideal in this semigroup and $u \in {\cal M}$ -- an idempotent. Then $u{\cal M}$ is the Ellis group of the flow $(G,S_{G,\ext}(M))$. Now, we adapt the notation from \cite[Chapter I.4]{Gl}. For any externally definable $G$-flow $(G,X)$ and a point $x_0 \in X$ such that $ux_0=x_0$, we define the Galois (or Ellis) group of $(X,x_0)$ as 

$${\cal G}(X,x_0)=\{ \alpha \in u{\cal M}: \alpha x_0=x_0\}. $$

Whenever there is a homomorphism $f$ between externally definable pointed flows $(X,x_0)$ and $(Y,y_0)$ (where $ux_0=x_0$ and $uy_0=y_0$), then ${\cal G}(X,x_0) \leq {\cal G}(Y,y_0)$.  Recall that a $G$-flow homomorphism is called proximal if any two points in each fiber are proximal. Then, $h$ is proximal if and only if ${\cal G}(X,x_0) = {\cal G}(Y,y_0)$.

We finish with a short discussion on equicontinuity. The essential notion here is uniformity. For the definition of uniformity and a concise exposition of basic material on this topic the reader is referred to the second appendix in \cite{Au}.

\begin{definition}
A $G$-flow $(G,X)$ is said to be equicontinuous if for any index $\alpha$ from the unique uniformity on $X$ there is an index $\beta$ such that $g\beta \in \alpha$ for all $g \in G$.
\end{definition}

Without giving the definition of uniformity in general, let us only say that in a compact Hausdorff space there is a unique uniformity generating the underlying topology, and it consists of all neighborhoods of the diagonal. In a compact Hausdorff group $K$, this uniformity has a basis consisting of the sets \{$(a,b) \in K \times K: (\exists U \in {\cal U})(a^{-1}b \in U)\}$, where ${\cal U}$ consists of all neighborhoods of the neutral element in $K$.

\section{Generalized Bohr compactification as a quotient of the Ellis group}\label{section 2}

As before, $G$ is a group $\emptyset$-definable in a model $M$ of an arbitrary theory, $N \succ M$ is $|M|^+$-saturated, and $\C \succ N$ is a monster model. We identify $S_{G,\ext}(M)$ with $S_{G,M}(N)$. Recall that the semigroup operation in $S_{G,\ext}(M)$ is denoted by $*$, but quite often we will just skip it. Moreover, ${\cal M}$ is a minimal subflow of $(G,S_{G,ext}(M))$ and $u \in {\cal M}$ is an idempotent. Thus, $u{\cal M}$ is the Ellis group of $S_{G,ext}(M)$.

In this section, we will generalize \cite[Chapter IX, Theorem 4.2]{Gl} on presenting the generalized Bohr compactification as a quotient of the Ellis group in our externally definable [or definable, assuming that all types in $S_G(M)$ are definable] context. The general approach via the so called $\tau$-topology is the same as in \cite[Chapter IX]{Gl}. However, there are some difficulties and the proof of the main theorem is different (although it uses various tricks and computations from \cite[Chapter IX]{Gl}). The reason is that we do not know whether the $G$-flow $2^{S_{G,\ext}(M)}$ of all non-empty closed subsets of $S_{G,\ext}(M)$ is externally  definable, and so we do not have a continuous on the left ``action'' of the universal externally definable $G$-ambit $(G,S_{G,\ext}(M),e)$ on this flow, and in consequence, we do not have all the nice properties of the circle operation considered in \cite[Chapter IX]{Gl}. We will be rather precise about what goes through and what does not.

Although we do not have the natural ``action'' of $S_{G,\ext}(M)$ on $2^{S_{G,\ext}(M)}$, we can take the statement in point (1) of \cite[Chapter IX, Lemma 1.1]{Gl} as the definition of $\circ$.

\begin{definition}
For $A \subseteq S_{G,\ext}(M)$ and $p\in S_{G,\ext}(M)$, $p \circ A$ is defined as the set of all points $x \in S_{G,\ext}(M)$ for which there exist nets $(x_i)$ in $A$ and $(g_i)$ in $G$ such that $\lim g_i=p$ and $\lim g_ix_i=x$.
\end{definition}

It is easy to check that $p\circ A$ is closed, $pA \subseteq p \circ A$ and $p\circ (q \circ A) \subseteq (pq)\circ A$ (but we do not know whether $p\circ (q \circ A) = (pq)\circ A$; in topological dynamics, it follows from the existence of the action of $\beta G$ on $2^{\beta G}$). Let us check for example that $p\circ (q \circ A) \subseteq (pq)\circ A$. 

Take $x \in p\circ (q \circ A)$. Then, there are nets $(g_i)$ and $(x_i)$ of elements of $G$ and $q \circ A$, respectively,  such that $\lim g_i=p$ and $\lim g_i x_i =x$. Next, for every $i$ there are nets $(h_{i,j})$ and $(y_{i,j})$ from $G$ and $A$, respectively, satisfying $\lim_j h_{i,j}=q$ and $\lim_j h_{i,j}y_{i,j}= x_i$. Take any neighborhood $U$ of $pq$ and $V$ of $x$. There is an index $i(U,V)$ such that $g_{i(U,V)}q \in U$ and $g_{i(U,V)}x_{i(U,V)} \in V$. So, there is an index $j(U,V)$ such that $g_{i(U,V)}h_{i(U,V),j(U,V)} \in U$ and $g_{i(U,V)}h_{i(U,V),j(U,V)}y_{i(U,V),j(U,V)} \in V$. All of this implies that $\lim_{U,V} g_{i(U,V)}h_{i(U,V),j(U,V)}= pq$ and $\lim_{U,V} (g_{i(U,V)}h_{i(U,V),j(U,V)})y_{i(U,V),j(U,V)}=x$.

\begin{definition}\label{Def: tau topology}
For $A \subseteq u{\cal M}$, define $\cl_\tau (A) = (u \circ A) \cap u{\cal M}$.
\end{definition}

Now, the proofs of 1.2-1.12 (except 1.12(2)) from \cite[Chapter IX]{Gl} go through (with some slight modifications) in our context. So we have all these results at our disposal. We also have a ``definable'' variant of Theorem 2.1 from there. 
%%%Krzys1: Added "externally".
However, a problem appears in Chapter IX.3, as for a $\tau$-closed subgroup $F$ of $u{\cal M}$ we do not have the canonical externally definable $G$-flow ${\cal A} (F):=\{ p \circ F : p \in {\cal M}\} \subseteq 2^{\cal M}$ whose Galois group at $u\circ F$ is $F$ (because we do not know whether $2^{\cal M}$ is an externally definable $G$-flow). Instead of this, we will use some other argument (based on \cite[Chapter IX, Proposition 4.1]{Gl} and some tricks from Chapter IX of \cite{Gl}). Let us note that Lemma 3.4 from Chapter IX also goes through in our context.

We will not repeat here all the necessary results from \cite{Gl}. Let us only recall a few basic facts. 
Namely, $\cl_\tau$ is a closure operator on subsets of $u{\cal M}$, and it induces the so-called $\tau$-topology on $u{\cal M}$. This topology is compact and $T_1$, and multiplication is continuous on each coordinate separately. 

\begin{definition}
$H(u{\cal M})$ is the intersection of the sets $\cl_\tau(V)$ with $V$ ranging over all $\tau$-neighborhoods of $u$ in the group $u{\cal M}$. 
\end{definition}
Then $H(u{\cal M})$ is a $\tau$-closed, normal subgroup of $u{\cal M}$, $u{\cal M}/H(u{\cal M})$ is a compact Hausdorff group, and for any $\tau$-closed subgroup $K$ of $u{\cal M}$, $u{\cal M}/K$ is a Hausdorff space if and only if $K \supseteq H(u{\cal M})$.

Our goal in this section is to show that $u{\cal M}/H(u{\cal M})$ is the generalized externally definable Bohr compactification of $G$.

%%%Krzys1: Added "externally" before "definable".
\begin{lemma}\label{sim}
Define $\sim$ on ${\cal M}$ by: 
$$p \sim q \iff p \in q \circ u{\cal M} \; \wedge \; q \in p \circ u{\cal M} \; \wedge \; up^{-1}q \in H(u{\cal M}).$$
Then $\sim$ is a closed, $G$-invariant equivalence relation on ${\cal M}$. So $(G, {\cal M}/\!\!\sim)$ is a minimal externally definable $G$-flow. Moreover, ${\cal G}({\cal M}/\!\!\sim,u/\!\!\sim)=H(u{\cal M})$. 
\end{lemma}

Note that whenever we compute the inverse $p^{-1}$ of an element $p \in {\cal M}$, we do this inside the group $v{\cal M}$, where $v$ is the unique idempotent such that $p \in v{\cal M}$.\\

\noindent
{\em Proof.} The fact that $\sim$ is a $G$-invariant equivalence relation uses the argument from the first paragraph of the proof of \cite[Chapter IX , Proposition 4.1]{Gl}. The fact that $\sim$ is closed is an exercise on limits of nets, and it also uses the argument from the second paragraph of that proof. External definability of $(G,{\cal M}/\sim)$ follows from the external definability of $(G,{\cal M})$. Let us show now that ${\cal G}({\cal M}/\!\!\sim,u/\!\!\sim)=H(u{\cal M})$.

Take $a \in u{\cal M}$ such that $a(u/\!\! \sim)=u/\!\! \sim$, i.e., $au/\!\! \sim =u/\!\! \sim$. Then $au \sim u$, so $u(u^{-1}(au)) \in H(u{\cal M})$, so $a \in H(u{\cal M})$. To prove the opposite inclusion, consider any $a \in H(u{\cal M})$. Then $uu^{-1}(au)=a \in H(u{\cal M})$. Moreover, $a \in u{\cal M}=u(u{\cal M}) \subseteq u \circ u{\cal M}$ and $ u \in a(u{\cal M}) \subseteq a \circ u{\cal M}$. All of this means that $u \sim au$, so $a(u/\!\! \sim)=u/\!\! \sim$, which completes the proof. \hfill $\square$

\begin{theorem}\label{generalized Bohr compactification as a quotient} 
${\cal M}/\!\! \sim$ is the universal externally definable compactification flow of $G$, and  the compact group $u{\cal M}/H(u{\cal M})$ is the generalized externally definable Bohr compactification of $G$.
\end{theorem}
{\em Proof.}
In this proof, instead of $H(u{\cal M})$ we will write $H$. We start from the following claim which is crucial.\\[3mm]
{\bf Claim} The compact group $u{\cal M}/H(u{\cal M})$ acts faithfully, jointly continuously and by automorphisms on the right on the flow ${\cal M}/\!\! \sim$ in the following way 
$$[p(u/\!\! \sim)] \bullet (aH)=p(a(u/\!\!\sim))=pa/\!\!\sim.$$
{\em Proof of the claim.} First, we check that $\bullet$ is well-defined, i.e., that it does not depend on the choice of $a$ in the coset $aH$ and that  it does not depend on the choice of $p$. The first statement follows from Lemma \ref{sim}, more precisely, from the fact that $H$ fixes $u/\!\!\sim$. For the second statement, consider any $p,q \in {\cal M}$ such that $p(u/\!\!\sim)=q(u/\!\!\sim)$. We need to show that $pa/\!\!\sim =qa/\!\!\sim$, i.e., $pa \sim qa$. Since $p \sim q$, we have that $up^{-1}q \in H$.
%$u(pa)^{-1}(qa)=ua^{-1}p^{-1}qa=uua^{-1}up^{-1}qua=u(ua)^{-1}(up^{-1}q)(ua) \in H$ 
Thus, $u(pa)^{-1}(qa)=a^{-1}p^{-1}qa=a^{-1}(up^{-1}q)a \in H$ (as $H$ is a normal subgroup of $u{\cal M}$ and $a \in u{\cal M}$). It remains to check that $pa \in qa \circ u{\cal M}$ (the symmetric condition follows in the same way). Since $p \sim q$, there are $g_i \in G$ and $x_i \in u {\cal M}$ such that $\lim g_i = q$ and $\lim g_ix_i = p$. Then, $pa=\lim (g_ia) (a^{-1}x_ia)$ and $\lim g_ia=qa$. Thus, one can find nets $(g_j')$ in $G$ and $(y_j)$ in $u{\cal M}$ such that $\lim g_j'=qa$ and $\lim g_j'(a^{-1}y_ja) =pa$, and hence $pa \in qa \circ u{\cal M}$.

Now, we check that for any $aH \in u{\cal M}/H$, $\bullet (aH)$ is an automorphism of the flow ${\cal M}/\!\!\sim$. $G$-invariance is clear: $[g(p(u/\!\!\sim))]\bullet (aH)=gpa/\!\!\sim = g([p(u/\!\!\sim)]\bullet (aH))$ for any $g \in G$ and $p \in {\cal M}$. To show continuity, consider any net $(p_i)$ in ${\cal M}$ such that $\lim p_i(u/\!\!\sim) = p(u/\!\!\sim)$ and suppose for a contradiction that $\lim p_ia/\!\!\sim \ne pa/\!\!\sim$. There exists a subnet $(q_j)$ of $(p_i)$ such that $\lim q_ja/\!\!\sim$ exists and is different from $pa/\!\!\sim$. There is a subnet $(s_k)$ of $(q_j)$ such that $\lim s_k$ exists in ${\cal M}$, and so $\lim s_k a \notin pa/\!\!\sim$. But $(s_k)$ is a subnet of $(p_i)$, so $\lim s_k \in p/\!\!\sim$, 
%%%Krzys1: I replaced "one can check that" by "by the previous paragraph".
hence, by the previous paragraph, $\lim s_ka \in pa/\!\!\sim$, a contradiction.

The fact that $\bullet$ is an action of $u{\cal M}/H(u{\cal M})$ on ${\cal M}/\!\! \sim$ is clear from the definition. It is faithful, because if $\bullet (aH)$ is a trivial automorphism, then $a \sim u$, so $a=uu^{-1}a \in H$.

We have  already checked that the action $\bullet$ is continuous in the first coordinate. In order to see that it is continuous as a 2-variable function, by the joint continuity theorem (see \cite[Chapter 4, Theorem 7]{Au}), it remains to show that it is continuous in the second coordinate. For this, it is enough to check that the function $f \colon u{\cal M} \to ({\cal M}/\!\!\sim)^{{\cal M}/\!\!\sim}$ given by the formula $f(a)(p/\!\!\sim)=pa/\!\!\sim$ is continuous. So, consider any net $(a_i)$ in $u{\cal M}$ such that $\tau$-$\lim a_i=a\in u{\cal M}$. Our goal is to show that $\lim pa_i/\!\!\sim =pa/\!\!\sim$ for any $p \in {\cal M}$. We will be done if we prove that for any subnet $(b_j)$ of $(a_i)$ for which $(pb_j)$ converges in the usual topology on ${\cal M}$ to some $r$ one has $r \sim pa$. Using Lemmas 1.4 and 1.5 from \cite[Chapter IX]{Gl}, the assumption that $\lim pb_j=r$ implies that $\tau$-$\lim upb_j=ur$, so, by \cite[Chapter IX, Proposition 1.7]{Gl}, $\tau$-$\lim ur^{-1}pb_j=\tau\mbox{-}\lim (ur^{-1})(upb_j)=ur^{-1}ur=ur^{-1}r=u$. On the other hand, $\tau$-$\lim ur^{-1}pb_j=ur^{-1}pa$. Therefore, since $u{\cal M}/H$ is Hausdorff, we conclude that $ur^{-1}paH=uH=H$, so 
%$$(*)\;\;\;\;\; 
\begin{equation}\label{equation *}
ur^{-1}pa \in H.
\end{equation}
Since $b_j, a \in u{\cal M}$, we easily get 
%$$(**)\;\;\;\;\; 
\begin{equation}\label{equation **}
r = \lim pb_j \in p \circ u{\cal M}=pa \circ u{\cal M}.
\end{equation}
Consider any open neighborhoods $U$ of $r$ and $V$ of $pa$. Choose $j(U,V)$ such that $pb_{j(U,V)} \in U$ (it exists as $\lim pb_j=r$). Choose $g_{j(U,V)} \in U \cap G$ so close to $pb_{j(U,V)}$ that $g_{j(U,V)} b_{j(U,V)}^{-1}a \in V$ (it exists as $pb_{j(U,V)}b_{j(U,V)}^{-1}a = pa$ and multiplication in ${\cal M}$ is continuous in the first coordinate). Then $\lim_{U,V} g_{j(U,V)}= r$ and $\lim _{U,V} g_{j(U,V)} b_{j(U,V)}^{-1}a=pa$, so
%$$(***)\;\;\;\;\; 
\begin{equation}\label{equation ***}
pa \in r \circ u{\cal M}.
\end{equation}
By (\ref{equation *}), (\ref{equation **}) and (\ref{equation ***}), we have $r \sim pa$, which was our goal, and the proof of the claim is completed. \hfill $\square$\\

Let $(G,\pi_{\extdef},x)$ be the universal minimal externally definable proximal $G$-flow (where $ux=x$), and let $\varphi \colon  ({\cal M},u) \to (\pi_{\extdef},x)$ be the epimorphism given by $\varphi(p)=px$ (recall that $px :=\lim g_ix$ for an arbitrary net $(g_i)$ in $G$ converging to $p$). 

We claim that $\varphi$ factors through the quotient map $i: ({\cal M},u) \to ({\cal M}/\!\!\sim, u/\!\!\sim)$. (Indeed, if $p \sim q$, then $p \in q \circ u{\cal M}$, so $p=\lim g_ix_i$ for some nets $(x_i)$ in  $u{\cal M}$ and $(g_i)$ in $G$ satisfying $\lim g_i=q$, and so $px=\lim g_i x_ix=\lim g_ix=qx$ (we use here the fact that ${\cal G}(\pi_{\extdef},x)=u{\cal M}$ as $(G,\pi_{\extdef})$ is proximal).) So, let $\bar \varphi \colon {\cal M}/\!\!\sim \to \pi_{\extdef}$ be the map such that $\bar \varphi i = \varphi$ (i.e. $\bar \varphi(p/\!\!\sim)=px$). 

Define 
$${\cal N}:=({\cal M}/\!\!\sim) / (u{\cal M}/H),$$ 
the space of orbits on ${\cal M}/\!\!\sim$ under the action $\bullet$ of $u{\cal M}/H$. Note that by the claim, ${\cal N}$ is a minimal externally definable $G$-flow.

Directly from the definition of ${\cal N}$ (and the claim), we see that $({\cal M}/\!\!\sim,u{\cal M}/H)$ is a group extension of ${\cal N}$. We will prove that:
\begin{enumerate}
\item ${\cal N} \cong \pi_{\extdef}$.
\item $\bar \varphi \colon ({\cal M}/\!\!\sim,u{\cal M}/H,u/\!\!\sim) \to (\pi_{\extdef},x)$ is the universal externally definable group extension of $\pi_{\extdef}$.
\end{enumerate}

From this it follows that ${\cal M}/\!\!\sim$ is the universal externally definable compactification flow of $G$ and that $u{\cal M}/H$ is the generalized externally definable Bohr compactification of $G$.

Note that $\bar \varphi$ factors through the quotient map $j: {\cal M}/\!\!\sim \to {\cal N}$. (Indeed, if $q \sim pa$ for some $a \in u \cal {M}$, then $qx=pax=px$). So, let $\bar {\bar \varphi} \colon {\cal N} \to \pi_{\extdef}$ be such that $\bar {\bar \varphi} j=\bar \varphi$.

Since ${\cal G}({\cal N},j(i(u)))=u{\cal M}={\cal G}(\pi_{\extdef},x)$, by \cite[Chapter I, Proposition 4.1(2)]{Gl}, we get that $\bar {\bar \varphi}$ is a proximal epimorphism, so, as $(G,\pi_{\extdef})$ is a proximal flow, ${\cal N}$ is also proximal. Since minimal proximal flows do not have endomorphisms different from the identity (see Lemma 4.1 in \cite[Chapter II]{Gl}) and $\pi_{\extdef}$ is the universal externally definable minimal proximal flow, we conclude that $\bar {\bar \varphi}$ is an isomorphism, so (i) is proved.

Since $j\colon ({\cal M}/\!\!\sim,u{\cal M}/H) \to {\cal N}$ is an externally  definable group extension and $\bar {\bar \varphi}$ is an isomorphism satisfying $\bar {\bar \varphi} j=\bar \varphi$, we get that $\bar \varphi \colon ({\cal M}/\!\!\sim,u{\cal M}/H,u/\!\!\sim) \to (\pi_{\extdef},x)$ is an externally  definable group extension. Finally,  the universality of this extension follows from the obvious counterpart of \cite[Chapter IX, Theorem 2.1(3)]{Gl}, the fact that group extensions are distal (i.e. any two points in each fiber are either equal or not proximal) and the equality ${\cal G}({\cal M}/\!\!\sim, u/\!\!\sim)=H$. \hfill $\square$

\section{Generalized Bohr compactification and connected components}\label{proof of main theorem 1}

The goal of this section is to prove Theorem \ref{main theorem 1}. But first we need to define the epimorphisms $\theta $ and $f$ discussed before this theorem.

We take the notation from Section \ref{section 2}. So we will be working in the general externally definable context. Recall that this applies to the definable case (assuming that all types in $S_G(M)$ are definable), and, in particular, to the case when predicates for all subsets of $G$ are in the language (so to classical topological dynamics). Theorem \ref{main theorem 1} is new even in this last situation. 

Take any set $A \subseteq M$.
Newelski \cite{Ne1} proved that the function $\theta:u{\cal M} \to {G^*}/{G^*}^{00}_A$ given by $f(\tp(a/N))=a{G^*}^{00}_A$ is a well-defined group epimorphism. One can easily refine this in the following way.

\begin{proposition}\label{definition of f}
The function $f:u{\cal M} \to {G^*}/{G^*}^{000}_A$ given by $f(\tp(a/N))=a{G^*}^{000}_A$ is a well-defined group epimorphism. Moreover, the composition of $f$ with the natural map from ${G^*}/{G^*}^{000}_A$ to ${G^*}/{G^*}^{00}_A$ is exactly $\theta$.
\end{proposition}
{\em Proof.} 
Consider $\hat{f}: S_{G,M}(N) \to G^*/{G^*}^{000}_A$ given by the same formula as $f$. First note that $\hat{f}$ is well-defined, because if $\tp(a/M)=\tp(b/M)$, then $b^{-1}a \in {G^*}^{000}_A$. 

Now, we check that $\hat{f}$ is a semigroup homomorphism. Take any $p,q \in S_{G,M}(N)$. By Fact \ref{description of *}, $p * q= \tp(ab/N)$, where $\tp(b/N)=q$ and $\tp(a/N,b)$ is the unique extension of $p$ to a complete type over $N,b$ which is finitely satisfiable in $M$. So, $\hat{f}(p)=a{G^*}^{000}_A$, $\hat{f}(q)=b{G^*}^{000}_A$ and $\hat{f}(p*q)=ab{G^*}^{000}_A$.

It remains to prove that $f$ is onto. Since each type in $S_G(M)$ can be extended to a type in $S_{G,M}(N)$, we easily get that $\hat{f}$ is onto. 
%Note that $\hat{f}(u)=e{G^*}^{000}_M$ (because $\hat{f}(u)=\hat{f}(u*u)=\hat{f}(u)^2$). 
We already know that $\hat{f}$ is a homomorphism. Thus, as ${\cal M}=S_{G,M}(N) * u$, we conclude that $\hat{f}|_{\cal M}$ is onto, and so $f=\hat{f}|_{u*{\cal M}}$ is also onto. \hfill $\square$\\

We have explained what the function $f$ in Theorem \ref{main theorem 1} is, and now we are ready to prove this theorem.\\

\noindent
{\em Proof of Theorem \ref{main theorem 1}.} 
(1) Let $\bar D \subseteq G^*/{G^*}^{000}_A$ be closed, i.e., $D:=\{ g \in G^*: g{G^*}^{000}_A \in \bar D\}$ is type-definable. The goal is to show that $f^{-1}[\bar D]$ is a $\tau$-closed subset of $u{\cal M}$.

Consider any $p \in \cl_{\tau}(f^{-1}[\bar D])$. By the definition of the $\tau$-topology, there are  $g_i \in G$ and $p_i \in f^{-1}[\bar D]$ such that $\lim_i g_i =u$ and $\lim_ig_ip_i =p$ (these limits are in the usual topology on $S_G(N)$). As $f^{-1}[\bar D]= \{ \tp(a/N): a \in D\} \cap u{\cal M}$ and $g_ip_i = \tp(g_ia_i/N)$ for $a_i \models p_i$, by compactness, we get that there are $a \models u$ and  $b \in D$ such that $ab \models p$. Then, since $f(u)=e{G^*}^{000}_A$,  we have $a \in {G^*}^{000}_A$. As $b \in D$ and $D$ is closed under multiplication by elements of ${G^*}^{000}_A$, we get $ab \in D$, so $p =\tp(ab/N)$ belongs to $f^{-1}[\bar D]$.\\[1mm]
(2) By Fact \ref{G000 as products}, we know that ${G^*}^{000}_A = \bigcup_{n \in \omega} F_n$, where $F_n$ is the $A$-type-definable set consisting of products of $n$ (equivalently, at most $n$) elements of the form $b^{-1}a$ where $(a,b)$ extends to an $A$-indiscernible sequence.
Put $\widetilde{F}_n:=\{ \tp(a/N): a \in F_n\}$, a closed subset of $S_G(N)$. Since $f(u)=e{G^*}^{000}_A$, we can find $n \in \omega$ such that $u \in \widetilde{F}_n$.

Let $\pi$ be the partial type over $A$ defining $F_{2n}$ and closed under conjunction.
Consider any $\varphi(x) \in \pi$. Let $$V:=[\neg \varphi(x)] \cap u{\cal M},$$ where  $[\neg \varphi(x)]$ is the closed subset of $S_G(N)$ consisting of types containing $\neg \varphi(x)$.\\[3mm]
{\bf Claim} \\
i) $u \notin \cl_{\tau}(V)$.\\
ii) $\cl_\tau (u{\cal M} \setminus \cl_\tau(V)) \subseteq \cl_\tau(u{\cal M} \setminus V) \subseteq \widetilde{F}_{3n}^\varphi$, where $\widetilde{F}_{3n}^\varphi := \{ \tp(ab/N): a \in F_n \wedge b \models \varphi(x)\}$.\\[3mm]
{\em Proof of the claim.}
i) %If $V$ is empty, there is nothing to do, so assume that $V \ne \emptyset$. 
Suppose for a contradiction that $u \in \cl_{\tau}(V)$. 
Arguing as in (1), we get that there are $a \models u$ and $b \models \neg \varphi(x)$ such that $ab \models u$. Then $a \in F_n$ and $ab \in F_n$, so $b \in F_{2n}$, and hence $b\models \varphi(x)$, a contradiction.\\[1mm]
ii) We need to check that  $\cl_\tau(u{\cal M} \setminus V) \subseteq \widetilde{F}_{3n}^\varphi$. Consider any $p \in \cl_\tau(u{\cal M} \setminus V)$. As before, there are $a \models u$ and $b \models \varphi(x)$ such that $ab \models p$. Then $a \in F_n$, so $\tp(ab/N) \in \widetilde{F}_{3n}^\varphi$. \hfill $\square$\\

Notice that $\bigcap_{\varphi(x) \in \pi}\widetilde{F}_{3n}^\varphi = \widetilde{F}_{3n}$. So, by the claim,
%%%Krzys1: I replaced $V$ by $U$ in the next formula, because $V$ was reserved for something else.
$$
\begin{array}{lll}
H(u{\cal M}) &= &\bigcap \left\{\cl_{\tau}(U): U \; \mbox{$\tau$-neighborhood of} \; u\right\} \subseteq \bigcap_{\varphi(x) \in \pi}\widetilde{F}_{3n}^\varphi \cap u{\cal M} = \widetilde{F}_{3n} \cap u{\cal M}\\ 
&\subseteq & f^{-1}(e{G^*}^{000}_M),
\end{array}$$
which finishes the proof of (2).\\[1mm]
(3) follows from (1) and (2).\\[1mm]
\indent
The fact that $u{\cal M}/H(u{\cal M})$ is the generalized externally definable Bohr compactification of $G$ was proved in Theorem \ref{generalized Bohr compactification as a quotient}. \hfill $\square$\\

Note that from Theorem \ref{main theorem 1} it follows immediately that the epimorphism $\theta$ from $u{\cal M}$ to ${G^*}/{G^*}^{00}_A$ is also continuous.

\section{The proof of Theorem \ref{main theorem 2}}

This section is devoted to the proof of Theorem \ref{main theorem 2}, which will consist of a series of lemmas.\\
%The context is still completely general (as in previous sections).

\noindent
{\em Proof of Theorem \ref{main theorem 2}.} Enlarging $M$ if necessary, we can assume that $A \subseteq M$, and we aim at proving the second (more precise) statement of the theorem. We take the notation as in Section \ref{section 2}. 

The situation is as follows.
By Theorem \ref{main theorem 1}, we have the natural continuous epimorphism  $\bar f: u{\cal M}/H(u{\cal M}) \to  G^*/{G^*}^{000}_A$ given by $\bar f(pH(u{\cal M}))=f(p)$, where $f: u{\cal M} \to G^*/{G^*}^{000}_A$ is defined by $f(\tp(a/N))=a{G^*}^{000}_A$. Next, $Y:=\ker(\bar f)=\bar f^{-1}(e{G^*}^{000}_A)$, and $\cl_\tau(Y)$ is the closure of $Y$ computed inside the compact topological group $u{\cal M}/H(u{\cal M})$. Thus, $\bar f$ induces an isomorphism between $\cl_\tau (Y)/Y$ and $\bar f[\cl_\tau (Y)]$. To prove the theorem, we need to show that $\bar f[\cl_\tau(Y)]={G^*}^{00}_A/{G^*}^{000}_A$. Of course, $\bar f[\cl_\tau(Y)] \leq {G^*}^{00}_A/{G^*}^{000}_A$, because by Theorem \ref{main theorem 1}, the preimage by $\bar f$ of the closed subgroup ${G^*}^{00}_A/{G^*}^{000}_A$ is closed, so it contains $\cl_\tau(Y)$. The rest of the proof is devoted to the proof of the opposite inclusion.

Let $P_u$ be the preimage of $Y$ by the quotient map $\pi$ from $u{\cal M}$ to $u{\cal M}/H(u{\cal M})$, and let $S=\cl_\tau(P_u)$. Then $\pi[S] = \cl_\tau (Y)$, as $\pi$ is continuous and closed. Hence, $\bar f[\cl_\tau(Y)] =f[S]$. Since a well-known result (see e.g. \cite[Proposition 3.5]{GiNe}) says that ${G^*}^{00}_A/{G^*}^{000}_A$ is the closure of the neutral element in $G^*/{G^*}^{000}_A$, we conclude that in order to finish the proof, we need to show that $f[S]$ is closed in $G^*/{G^*}^{000}_A$.

All the sets $Y$, $S$ and $P_u$ are defined in terms of the idempotent $u \in {\cal M}$ chosen at the beginning. For any other idempotent $v \in {\cal M}$ we have the analogous sets, in particular, we have the set $P_v \subseteq v{\cal M}$.

Let $\hat{f} : S_{G,M}(N) \to G^*/{G^*}^{000}_A$ be the semigroup homomorphism extending $f$ given by $\hat{f}(\tp(a/N))=a{G^*}^{000}_A$ (see the proof of Proposition \ref{definition of f}). By $J$ we will denote the set of all idempotents in ${\cal M}$. By $S_{{G^*}^{000}_A,M}(N)$ we will, of course, denote the subset of $S_{G,M}(N)$ consisting of types of elements from ${G^*}^{000}_A$.

Below, sometimes we use $*$ to denote the semigroup operation in $S_{G,M}(N)$, but sometimes we do not use any symbol in such a situation.

%By the definition of $\tau$-topology, we easily get that $S=u(u \circ P)$, and so $f[S]=\hat{f}[u(u\circ P)]=\hat{f}[u \circ P]$. So the goal is to show that $\hat{f}^{-1}[\hat{f}[u \circ P]]$ is closed in the usual topology on $S_G(M)$. I did some computations around this, but could not get the conclusion. For example, I computed that $\hat{f}^{-1}[\hat{f}[u \circ P]]=\{q \in S_G(M): q * u \in S_{{G^*}^{000}_M}(M)*(u \circ P)\}$ (where $*$ is the semigroup operation). Hence, since $*$ is continuous on the first coordinate, our desired conclusion is equivalent to $S_{{G^*}^{000}_M}(M)*(u \circ P)$ being closed in $S_G(M)$. Is it true? One can also show that $\hat{f}^{-1}[\hat{f}[u \circ P]] \cap {\cal M} = J * (u \circ P) * S_{{G^*}^{000}_M}(M)*u$, where $J$ is the set of idempotents in ${\cal M}$. So if our conclusion holds, the last set should be closed. Is it true?\\

\begin{lemma}\label{lemat 1} For any idempotent $v \in J$ one has $P_v=\ker(\hat{f}) \cap v{\cal M} =vS_{{G^*}^{000}_A,M}(N)v$, and for any other idempotent $w \in J$ we have $P_v=vP_w$. Also, $P_v$ is a subgroup of $v{\cal M}$ for $v \in J$.
\end{lemma}
{\em Proof.} 
The first equality is obvious from the definition of $P_v$. The rest follows easily from the fact that $\hat{f}$ is a semigroup homomorphism mapping both $J$ and $S_{{G^*}^{000}_A,M}(N)$ to the neutral element of $G^*/{G^*}^{000}_A$. 

Let us check for example that $P_v=vP_w$. Since $\hat{f}[vP_w]=\hat{f}(v)\hat{f}[P_w]=\{ e{G^*}^{000}_A\}$, we get $P_v\supseteq vP_w$. For the opposite inclusion, take any $p \in P_v$. Then $p=vq$ for some $q \in {\cal M}$. Since $\hat{f}$ is a homomorphism and $v,w,p \in \ker (\hat{f})$, we get that $q,wq \in \ker{\hat{f}}$, so $wq \in P_w$ and $p=v(wq) \in vP_w$. \hfill $\square$

\begin{remark}
For any $D \subseteq G^*/{G^*}^{000}_A$ one has that $D$ is closed iff $\hat{f}^{-1}[D]$ is closed.
\end{remark}
{\em Proof.} %By the definition of the logic topology, the statement that $D$ is closed means that $E:=\{ \tp(a/N) : a{G^*}^{000}_A \in D\}$ is closed in $S_G(N)$, whereas the statement on the right hand side says that $E \cap S_{G,M}(N)$ is closed (so it is literally weaker). To resolve this problem, let us define a topology $T$ on $G^*/{G^*}^{000}_A$ by saying that $D$ is closed iff $\hat{f}^{-1}[D]$ is closed (in $S_{G,M}(N)$). Since $T$ is stronger than the logic topology, $T$ must be Hausdorff. In order to finish the proof, it is enough to see that $T$ coincides with the logic topology, and for that it is enough to show that $T$ is compact. But this follows easily from the fact that $\hat{f}$ is surjective and $S_{G,M}(N)$ is compact.\hfill $\square$
Since $A\subseteq M$, $D$ is closed iff $E:=\{ \tp(a/M) : a{G^{*}}^{000}_A \in D\}$ is closed in $S_G(M)$. The restriction map $\alpha: S_{G,M}(N) \to S_G(M)$ is continuous and closed. Moreover, $\hat{f}^{-1}[D] = \alpha^{-1}[E]$ and $\alpha[ \hat{f}^{-1}[D]]=E$. Thus, $D$ is closed iff $E$ is closed iff $\hat{f}^{-1}[D]$ is closed. \hfill $\square$

\begin{remark}\label{closeness  in M}
For any $D \subseteq G^*/{G^*}^{000}_A$ one has that $D$ is closed iff $\hat{f}^{-1}[D] \cap {\cal M}$ is closed in ${\cal M}$. In particular, in our context,  in order to prove that $f[S]$ is closed, it is enough to show that $\hat{f}^{-1}[f[S]] \cap {\cal M}$ is closed in ${\cal M}$.
\end{remark}
{\em Proof.} 
Note that $\hat{f}^{-1}[D] \cap {\cal M}=\hat{f}^{-1}[D]*u$. Moreover, $*$ is continuous in the left coordinate and $S_{G,M}(N)$ is compact Hausdorff. Thus, $\hat{f}^{-1}[D] \cap {\cal M}$ is closed iff $\hat{f}^{-1}[D]*u$ is closed iff $\hat{f}^{-1}[D]$ is closed (in $S_{G,M}(N)$) iff $D$ is closed (where the last equivalence follows from the previous remark). \hfill $\square$\\

So, in order to finish the proof of the theorem, it remains to show that $\hat{f}^{-1}[f[S]] \cap {\cal M}$ is closed in ${\cal M}$.

\begin{lemma}\label{lemat 2}
$\hat{f}^{-1}[f[S]] \cap {\cal M} = J * (u \circ P_u) * S_{{G^*}^{000}_A,M}(N)*u=J*(u \circ P_u)$.
\end{lemma}
{\em Proof.}
The second equality follows easily from the continuity of $*$ on the left and from the observation that $P_u*S_{{G^*}^{000}_A,M}(N)*u=P_u$ which is an immediate consequence of  Lemma \ref{lemat 1}. 

So, we concentrate on the first equality. 
The inclusion $(\supseteq)$ follows from the observation that $f[S]=\hat{f}[u \circ P_u]$ (being an obvious consequence of the fact that $S=u(u \circ P_u)$ noted in \cite[Chapter IX, Lemma 1.3]{Gl}) and from the fact that $\hat{f}$ is a semigroup homomorphism which maps all idempotents and all elements from $S_{{G^*}^{000}_A,M}(N)$ to the neutral element in $G^*/{G^*}^{000}_A$.

$(\subseteq)$ Take any $q \in \hat{f}^{-1}[f[S]] \cap {\cal M}$. Then $\hat{f}(q)=f(s)$ for some $s \in S=u(u \circ P_u)$. There is $v \in J$ such that $q \in v{\cal M}$. Then $\hat{f}(q)=\hat{f}(vs)$. So $\hat{f}((vs)^{-1} q) = e{G^*}^{000}_A$, i.e., $(vs)^{-1} q \in S_{{G^*}^{000}_A,M}(N) \cap {\cal M} = S_{{G^*}^{000}_A,M}(N)*u$. Therefore, $vq \in vs *S_{{G^*}^{000}_A,M}(N)*u$, but $vq=q$, so $q \in vs *S_{{G^*}^{000}_A,M}(M)*u \subseteq J * (u \circ P_u) * S_{{G^*}^{000}_A,M}(N)*u$. \hfill $\square$

\begin{lemma}\label{lemat 3} Let $v \in J$.\\
i) $v \circ P_u= v \circ P_v$ and $u \circ P_v= u \circ P_u$.\\
ii) $v(u \circ P_u) \subseteq v \circ P_v$ and $u(v \circ P_v) \subseteq u \circ P_u$.\\
iii) $v(u \circ P_u)=v(v \circ P_v)$ and $u( v\circ P_v)=u(u \circ P_u)$.
\end{lemma}
{\em Proof.} 
In all three points, it is clearly enough to show only one of the two symmetric statements.\\[1mm]
Point (i) is a consequence of the following two sequences of inclusions which follow from Lemma \ref{lemat 1} and basic properties of $\circ$.
$$v \circ P_v = v \circ (vP_u) \subseteq v \circ(v\circ P_u) \subseteq (vv)\circ P_u=v \circ P_u$$ and
$$ v \circ P_u=v \circ (uvP_u)=v \circ (uP_v) \subseteq v \circ(u\circ P_v) \subseteq (vu) \circ P_v=v \circ P_v.$$
ii) By (i), we have 
$$v(u \circ P_u) \subseteq v \circ (u \circ P_u) \subseteq (vu) \circ P_u=v \circ P_u=v \circ P_v.$$
iii) By (ii), we get $v(u \circ P_u) =v(v(u \circ P_u)) \subseteq v(v \circ P_v)$ and similarly  $u( v\circ P_v)\subseteq u(u \circ P_u)$. It follows that $u(u \circ P_u)=u(v(u \circ P_u)) \subseteq uv(v \circ P_v) \subseteq u(u \circ P_u)$. So the second inclusion must be the equality, i.e., $u(v \circ P_v) = u(u \circ P_u)$. 
\hfill $\square$\\

The next lemma is a corollary of the above considerations.

\begin{lemma}\label{lemat 4}
$\hat{f}^{-1}[f[S]] \cap {\cal M}=  \bigcup_{v \in J} v \circ P_v=\bigcup_{v \in J} v \circ P_u$.
\end{lemma}
{\em Proof.} 
The second equality follows from Lemma \ref{lemat 3}(i). For the inclusion $(\subseteq)$ in the first equality, note that by Lemmas \ref{lemat 2} and \ref{lemat 3}(ii), 
$$\hat{f}^{-1}[f[S]] \cap {\cal M}=\bigcup_{v \in J} v (u\circ P_u) \subseteq \bigcup_{v \in J} v \circ P_v.$$
For $(\supseteq)$ notice that by Lemma \ref{lemat 3}(iii), we have $\hat{f}[v \circ P_v] =\hat{f}[v(v \circ P_v)]=\hat{f}[v(u \circ P_u)]=\hat{f}[u \circ P_u]=f[S]$, and so $v \circ P_v \subseteq \hat{f}^{-1}[f[S]] \cap {\cal M}$. \hfill $\square$

\begin{lemma}\label{lemat 5}
$\cl (J) \subseteq \ker(\hat{f}) \cap {\cal M}$. Equivalently, $\cl (J) \subseteq \bigcup_{v \in J} P_v = \bigcup_{v \in J} vP_u$.
\end{lemma}
{\em Proof.} 
The second part is equivalent to the first one, because $\bigcup_{v \in J} v{\cal M}={\cal M}$ and $P_v = \ker(\hat{f}) \cap v{\cal M}$ for $v \in J$. To show the first part, consider any $p \in \cl(J)$. 

For an arbitrary $v \in J$, we have that $v * v=v$. Thus, by Fact \ref{description of *}, for any $b \models v$ and $a$ realizing the unique extension of $v$ to a complete type over $N,b$ finitely satisfiable in $M$, we have that $ab \models v$. So, denoting by $F(x)$ the type over $M$ saying that $x=yz^{-1}$ for some $y \equiv_M z$, we get that $a=abb^{-1} \models F$. So, for any $v \in J$ and $c \models v$ we have that $c \models F$.

Since in every open neighborhood of $p$ there is an element from $J$, we have that for every $\varphi(x) \in p$ there is $c \models \varphi(x) \wedge F(x)$. By compactness, there is $c \models p(x) \wedge F(x)$. So $c \in {G^*}^{000}_A$ and $p \in \hat{f}^{-1}(e{G^*}^{000}_A) \cap {\cal M}$, and we are done. \hfill $\square$

\begin{lemma}\label{lemat 6}
$\bigcup_{v \in J} v \circ P_u$ is closed.
\end{lemma}
{\em Proof.} 
Consider any $q \in \cl(\bigcup_{v \in J} v \circ P_u)$. Then there are nets $(v_i)$ in $J$ and $(x_i)$ in ${\cal M}$ such that $x_i \in v_i \circ P_u$ and $(x_i)$ converges to $q$. We can find subnets $(w_j)$ of $(v_i)$ and $(y_j)$ of $(x_i)$ for which $(w_j)$ converges to some $r \in \cl(J)$ (and still $y_j \in w_j \circ P_u$). Thus, we can find nets $(g_k)$ in $G$ and $(p_k)$ in $P_u$ such that $(g_k)$ converges to $r$ and $(g_kp_k)$ converges to $q$. So we have proved that $q \in r \circ P_u$.

%%%Krzys1: I replaced $v$ by $w$ to distinguish it from $v$ appearing under sums.
By Lemma \ref{lemat 5}, we get that $r=wp$ for some $w \in J$ and $p \in P_u$. Since $P_u$ is a subgroup of $u{\cal M}$, we also have that $p^{-1}P_u=P_u$. Thus, we conclude that
$$q \in r \circ P_u = r \circ(p^{-1}P_u) \subseteq r \circ(p^{-1} \circ P_u) \subseteq (rp^{-1})\circ P_u=w \circ P_u \subseteq \bigcup_{v \in J} v \circ P_u,$$
which completes the proof of the lemma. \hfill $\square$\\

By Lemmas \ref{lemat 4} and \ref{lemat 6} together with the conclusion right after Remark \ref{closeness in M}, the proof of Theorem \ref{main theorem 2} has been completed. \hfill $\square$\\

Note that if one is not interested in proving that  ${G^*}^{00}_A/{G^*}^{000}_A$ is the quotient of a compact Hausdorff group by a dense subgroup but just by any subgroup, then such a weaker statement can be easily deduced from Theorem \ref{main theorem 1}. Namely, by this theorem, $Z:=\bar f^{-1}[{G^*}^{00}_A/{G^*}^{000}_A]$ is closed (so compact Hausdorff) in $u{\cal M}/H(u{\cal M})$, and clearly $\bar{f}$ induces an isomorphism between $Z/Y$ and ${G^*}^{00}_A/{G^*}^{000}_A$, where $Y=\ker (\bar f)$. In a similar fashion, we get

\begin{corollary}
For an arbitrary group $G$ which is $\emptyset$-definable in an arbitrary model $M$ and for any set of parameters $A$ the group $G^*/{G^*}^{000}_A$ is isomorphic to a quotient of a compact Hausdorff group. In the above notation, it will be isomorphic to the quotient of $u{\cal M}/H(u{\cal M})$ by $\ker({\bar f})$.
\end{corollary}

We finish this section proving Corollary \ref{easy corollary}.\\

\noindent
{\em Proof of Corollary \ref{easy corollary}.} 
(ii) follows from (i) using the fact that $\bar f$ is naturally induced by $f$. To see (i), notice that $Y:=\ker (\bar f)$ is a singleton, so it is closed, so $\cl_\tau(Y)=Y$, and, by Theorem \ref{main theorem 2}, $\{e{G^*}^{000}_A\}=\bar f [Y]= \bar f [\cl_\tau(Y)]={G^*}^{00}_A/{G^*}^{000}_A$, which implies that ${G^*}^{000}_A={G^*}^{00}_A$. \hfill $\square$

\section{Strongly amenable case}\label{strongly amenable case}

We will prove Corollary \ref{main corollary 3}, and then we will give some discussion around the proof.\\

\noindent
{\em Proof of Corollary \ref{main corollary 3}.} 
Recall that we assume that all types in $S_G(M)$ are definable, so $S_{G,\ext}(M)=S_G(M)$, and so ${\cal M}$ is a minimal left ideal in $S_G(M)$ and $u \in {\cal M}$ is an idempotent. Let $\pi \circ \bar f:u{\cal M}/H(u{\cal M}) \to G^*/{G^*}^{00}_M$ be denoted by $\zeta$. It is an epimorphism and the goal is to prove that it is an isomorphism.

It is enough to show that there is a definable (in the sense of Definition \ref{definable functions}) homomorphism $\eta : G \to u{\cal M}/H(u{\cal M})$ with dense image and such that $(\zeta \circ \eta)(g)=g{G^*}^{00}_M$ for all $g \in G$. Indeed, since $G^*/{G^*}^{00}_M$ is the definable Bohr compactification of $G$, universality yields a continuous homomorphism $h: G^*/{G^*}^{00}_M \to u{\cal M}/H(u{\cal M})$ such that $h(g{G^*}^{00}_M)=\eta(g)$ for all $g \in G$; then $(\zeta \circ h) (g{G^*}^{00}_M)=g{G^*}^{00}_M$, so $\zeta \circ h=\id$, so $\zeta$ is injective (as $h$ is onto by the denseness of $\im(\eta)$), so we are done.

Let $(G,X,K,x_0)$ be the universal definable compactification flow of $G$ (where $ux_0=x_0$). So, by definable strong amenability, it is the universal group extension of the trivial $G$-flow.
By Theorem \ref{generalized Bohr compactification as a quotient} and its proof, we know that $K= u{\cal M}/H(u{\cal M})$, and it acts on the right on $X$ by $(px_0)k=p(kx_0)$ for $p \in {\cal M}$. Since $(G,X,K)$ is a group extension of the trivial $G$-flow, we have that $X=x_0K$, and for any $g \in G$, $x \in X$ and $k \in K$ one has $g(xk)=(gx)k$. Then $K$ is the group of all automorphisms of the $G$-flow $X$. Also, the function $k \mapsto x_0k$ is a homeomorphism between $K$ and $X$.

%Let $(G,X,K,x_0)$ be the universal group extension of the trivial $G$-flow (in particular, $ux_0=x_0$). By the definition of a group extension, there is a compact topological group $K$ acting continuously on $X$ on the right such that $X=x_0K$ and for any $g \in G$, $x \in X$ and $k \in K$ one has $g(xk)=(gx)k$. Then $K$ is the group of all automorphisms of the $G$-flow $X$.

%By Theorems 2.1 and 4.2 (and their proofs) from Chapter IX of Glasner's book, we know that the function $\widetilde{\kappa} : u{\cal M}/H(u{\cal M}) \to K$ given by the formula $(px_0)\widetilde{\kappa}(\beta)=p(\beta x_0)$ (for $p \in {\cal M}$)  is a well-defined topological isomorphism.

%Since $(gx)k=g(xk)$ for all $g \in G$, $x \in X$ and $k \in K$, 
One easily gets that $(G,X)$ is equicontinuous. Thus, by \cite[Chapter I, Theorem 3.3]{Gl}, the Ellis semigroup $E(X)$ is a compact topological group consisting of homeomorphisms. By the definition of $E(X)$, we also have a natural homomorphism $\chi : G \to E(X)$ (namely, $\chi(g)$ is the translation by $g$) with dense image. Since $E(X)$ is a subflow of the product $G$-flow $X^X$, by Remark \ref{definability of products}, we easily get that $\chi$ is definable.

Define a function $\alpha: E(X) \to X^X$ by
$$(px_0)\alpha(\beta)=p(\beta x_0)$$
for $p \in E(X)$. To see that $\alpha$ is well-defined, notice that whenever $px_0=qx_0$, then $p(x_0k)=(px_0)k=(qx_0)k=q(x_0k)$ for all $k \in K$, and hence $p=q$. This computation also shows that $\alpha$ is injective.\\[3mm]
{\bf Claim} $\alpha$ is a topological isomorphism from $E(X)$ onto $K$.\\[3mm]
{\em Proof of the claim.}
First, let us check that for every $\beta \in E(X)$, $\alpha(\beta)$ is an automorphism of $(G,X)$, and so $\alpha(\beta) \in K$. Of course, for $ g \in G$ and $p \in E(X)$ we have $(gpx_0)\alpha(\beta)=gp(\beta x_0)=g(p\beta x_0)=g((px_0)\alpha(\beta))$. Now, we show continuity of $\alpha(\beta)$. Assume $\lim_i p_ix_0=px_0$. Then $\lim_i(p_ix_0)k=(px_0)k$, and so $\lim_i p_i(x_0k)=p (x_0k)$, for all $k \in K$. But $x_0K=X$, so $\lim_i p_i(\beta x_0) =p(\beta x_0)$, i.e. $\lim_i (p_i x_0)\alpha(\beta)=(px_0) \alpha(\beta)$. Thus, $\alpha(\beta)$ is continuous.

The fact that $\alpha:E(X) \to K$ is a group monomorphism is obvious.  To see that $\alpha$ is onto, consider any $k \in K$. Then, since $(G,X)$ is a minimal flow, there is $\beta \in E(X)$ such that $x_0\alpha(\beta)=\beta x_0=x_0k$, so $\alpha(\beta)=k$.

It remains to check that $\alpha$ is continuous. For this, we need to show that for any $p \in E(X)$ and $U$ an open subset of $X$, the collection $B$ of all $\beta \in E(X)$ satisfying $(px_0)\alpha(\beta) \in U$ is open in $E(X)$.  But since $E(X)$ is group, we get that $\beta \in B$ iff $p\beta x_0 \in U$ iff $\beta x_0 \in p^{-1}U$. As $E(X)$ consists of homeomorphisms of $X$, $p^{-1}U$ is open in $X$, and so $B$ is open in $E(X)$. \hfill $\square$\\

Now, define $\eta: G \to u{\cal M}/H(u{\cal M})$ as $\alpha \circ \chi$. By the above observations that $\alpha$ is a topological isomorphism and $\chi$ is a definable homomorphism with dense image, we conclude that $\eta$ is a definable (group) homomorphism with dense image. To finish the proof of the corollary, we need to show that  $(\zeta \circ \eta)(g)=g{G^*}^{00}_M$. For this, it is enough to see that $\eta(g)=(u g u) H(u{\cal M})$ (because then, taking $b\models u$ and $a$ realizing the unique coheir of $u$ over $M,b$, we have $ugu=\tp(agb/M)$ and $a,b \in {G^*}^{00}_M$, so $\zeta(uguH(u{\cal M}))=agb {G^*}^{00}_M=g{G^*}^{00}_M$). 

Since $ux_0=x_0$, we have that $u(x_0k)=(ux_0)k=x_0k$ for all $k \in K$, and so $u$ acts trivially on $X$.  Thus, $x_0\alpha(\chi(g))=\chi(g)x_0=gx_0=ugux_0=x_0(uguH(u{\cal M}))$. So $\eta(g)=\alpha(\chi(g))=uguH(u{\cal M})$. \hfill $\square$\\

%%%Krzys1: I changed the next paragraph. Before I wrote that "it is also clear that in general the image of $\eta$ is dense", but now I cannot see this. So I removed this sentence, and I weakened the rest of the paragraph. It is interesting which property "\eta is a homomorphism" or "the image of \eta is dense" or "both of them" fails in general. Now, it is only noted that one of them fails. 
Note that we obtained an explicit formula for the function $\eta:  G \to u{\cal M}/H(u{\cal M})$ in the proof above, namely $\eta(g)=uguH(u{\cal M})$. As was easily checked, even without using strong amenability, this function satisfies $(\zeta \circ \eta) (g)= g{G^*}^{00}_M$. 
%it is also clear that in general the image of $\eta$ is dense. 
The whole difficulty was to get that $\eta$ is a definable homomorphism with dense image, and for that we used strong amenability. It is worth to emphasise that in general $\eta$ will not be a homomorphism with dense image. To see that, it is enough to take any group $G$, with predicates for all subsets added to the language, such that ${G^*}^{00}_M \ne {G^*}^{000}_M$ (e.g. as $G$ one can take the free group in at least 2 generators \cite{GiKr}). Then clearly $\zeta$ is not injective, but if $\eta$ was a homomorphism with dense image, then the second paragraph of the above proof would imply that $\zeta$ is injective, a contradiction.

Since all nilpotent groups are strongly amenable, Corollary \ref{main corollary 4} follows immediately from Corollary \ref{main corollary 3}. The solvable case remains open. The following general and related question looks interesting.

%%%Krzys1: Do we really want to work without the assumption that all types in $S_G(M)$ are definable in point (i) of the next question? Then even strongly amenable case is not clear. So I added this assumption, and in consequence removed the reference to \cite{ChPiSi} in the comment below. Is it OK with you?
\begin{question} 
%(i) Suppose $G$ is a group definable in a structure $M$ and $G$ is definably amenable (in the sense of Definition \ref{definition of strong amenability}). Is it the case that ${G^*}^{00}_M={G^*}^{000}_M$?\\
(i) Let $G$ be a group definable in a structure $M$. 
%Assume all types in $S_G(M)$ are definable and 
Suppose $G$ is definably amenable (in the sense of Definition \ref{definition of strong amenability}). Is it the case that ${G^*}^{00}_M={G^*}^{000}_M$?\\
(ii)  The special case of (i) where all subsets of $G$ are definable in $M$: Namely, if the discrete group $G$ is amenable, is it the case that the Bohr compactification of $G$ coincides with the ``new'' invariant 
$G^{*}/{G^*}^{000}_M$?
\end{question}

%%%Krzys: I extended the comment below.
%Note that (i) has a positive answer when $\Th(M)$ has NIP. This was proved in \cite{HrPi}, but it also follows from Corollary \ref{main corollary 3} and Theorem \ref{main theorem 5} (proved in the next section). 

It would be natural to make a definability of types assumption in (i), but
the question makes sense without it, and a positive answer was given in \cite{HrPi}
when $T$ has NIP.
But the NIP case can also be deduced
from Corollary \ref{main corollary 3} and Theorem \ref{main theorem 5} (proved in the next section), bearing in
mind the invariance of the relevant connected components when passing to
$M^{\ext}$ which was proved in \cite{ChPiSi} (see also the beginning of the next section for a short discussion on $M^{\ext}$).

%together with the results from \cite{ChPiSi} saying that in the NIP environment connected components do not change when passing to the expansion by predicates for all externally definable subsets of $M$.

%\begin{question}
%Let $G$ be a solvable group definable in a structure $M$ is which predicates for all subsets of $G$ are in the language. Is it true ${G^*}^{00}_M={G^*}^{000}_M$?
%\end{question}

\section{Strong amenability and amenability under NIP}

The main point of this section is to give a proof of Theorem \ref{main theorem 5}, but we also take the opportunity to recall some issues on definable measures.

See earlier papers such as \cite{HrPi} for the definition of the property NIP, but if it helps the reader, it means that any definable family of definable sets has finite Vapnik-Chervonenkis dimension. 

With any first order structure $M$, we can associate the structure $M^{\ext}$ with the same universe as $M$ and predicates for all externally definable subsets of $M$. Recall that Shelah proved that if $M$ is a model of a NIP theory, then $\Th(M^{ext})$ is also NIP, it has quantifier elimination, and all complete types over $M^{\ext}$  in this theory are definable. Therefore, if $G$ is a group definable in a model $M$ of a NIP theory, then externally definable strong amenability of $G$ is equivalent to definable strong amenability of $G$ treated as a group definable in $M^{\ext}$. On the other hand, by \cite[Theorem 3.17]{ChPiSi}, we know that in a NIP theory definable amenability is preserved under passing to $M^{\ext}$ and coincides with externally definable amenability. That is why, replacing $M$ by $M^{\ext}$, we can assume that all types in $S_G(M)$ are definable and we can rewrite Theorem \ref{main theorem 5} in the following equivalent form.

\begin{theorem}
Let $G$ be a group definable in a model $M$ of NIP theory. Assume that all types in $S_G(M)$ are definable. Then $G$ is definably strongly amenable if and only if it is definably amenable.
\end{theorem}

From now on, let $G$ be a group definable in a first order structure $M$, and assume that all types in $S_G(M)$ are definable. Let $T=\Th(M)$.

First, we will prove the implication $(\leftarrow)$ by showing the following:

\begin{proposition}\label{amenability implies strong amenability}
Assume that $T$ is NIP and $G$ is definably amenable. Suppose $P$ is a definable minimal proximal $G$-flow. Then $P$ is trivial, i.e., it is a single point. 
\end{proposition}

Then, after a discussion on the (compact) space ${\cal M}(S_G(M))$ of Borel probability measures on $S_G(M)$ and on definability of measures, we will obtain the implication $(\rightarrow)$ by showing the following two facts.

\begin{proposition}\label{definability of the flow of measures}
Suppose $T$ is $NIP$. Then ${\cal M}(S_{G}(M))$ is a definable $G$-flow. 
\end{proposition} 

\begin{proposition}\label{strong amenability implies amenability} Suppose that the natural action of $G$ on ${\cal M}(S_{G}(M))$ is definable (i.e. ${\cal M}(S_{G}(M))$ is a definable $G$-flow). Assume $G$ is definably strongly amenable. Then $G$ is definably amenable.
\end{proposition}
{\em Proof of Proposition \ref{amenability implies strong amenability}.}
We will use recent work of Chernikov and Simon \cite{ChSi}, and we also refer to Glasner's book \cite{Gl}.  

First we recall the Baire-generic compact domination theorem from \cite{ChSi}. Since $G$ is definably amenable, there is a strongly $f$-generic global type $q \in S_G(\C)$ (according to the terminology from \cite{ChSi}) which was proved in \cite{HrPi}. Let $C$ be the closure of the orbit $G^*q$ in $S_G(\C)$, and let $\pi: S_G(\C) \to G^*/{G^*}^{00}$ be the obvious continuous map. By $\pi |_C: C \to G^*/{G^*}^{00}$ we denote the restriction of $\pi$ to $C$. For a clopen subset $X$ of $S_G(\C)$ define
%
%\begin{equation}\label{Baire-generic}
$$E_X:= \{ h \in G^*/{G^*}^{00} : (\pi|_{C})^{-1}(h) \cap X \ne \emptyset \;\; \mbox{and}\;\;(\pi|_{C})^{-1}(h) \cap X^c \ne \emptyset\}.$$
%\end{equation}
%
Theorem 5.2 from \cite{ChSi} tells us that $E_X$ is closed and has empty interior. Now, we want to deduce a counterpart of this working in our non-saturated model $M$. Let $r : S_G(\C) \to S_G(M)$ be the restriction map. Then $r[C]$ is a $G$-subflow of the flow $(G,S_G(M))$. Choose a minimal subflow $\mathcal{M}$ of $S_G(M)$ which is contained in $r[C]$.  Let now $\theta:{\cal M} \to G^*/{G^*}^{00}$ be the canonical, continuous, surjective homomorphism (it was defined on $u{\cal M}$ in Section \ref{proof of main theorem 1}, but the same definition works on ${\cal M}$). Then $\theta \circ (r|_{r^{-1}(\mathcal{M})})=\pi|_{r^{-1}(\mathcal{M})}$. Thus, by the fact that $E_X$ has empty interior for any clopen $X \subseteq S_G(\C)$, we easily get the following claim. (Closedness follows from the continuity of $\theta$.)\\

 \noindent
 {\bf Claim 1.} For every clopen subset $X$ of $\cal M$ (i.e. given by a formula over $M$) the set $\{g\in G^*/{G^*}^{00}: \theta^{-1}(g)$ meets both $X$ and $X^{c}\}$ is closed and has empty interior.\\

%%%Krzys: I do not have Chernikow and Simons's paper, so I did not check the places in which we refer to this paper. Is the information on $f$-generics relevant here? Maybe we could remove it.
Let $x_{0}\in P$. 
%and let  ${\cal M}$ be any minimal subflow of $S_{G}(M)$, which we know to consist of some $f$-generic types (in the sense of \cite{ChSi}).  Let $u$ be an idempotent in ${\cal M}$. 
%
%%%Krzys: I modified slightly this part, referring to the section with preliminaries.
By Fact \ref{universal definable ambits} and the comments below it, we know that $S_G(M)$ acts on $P$ in the sense that $ex=x$ and $(p *q)x=p(qx)$ for all $x \in X$ and $p,q \in S_G(M)$. Let $\phi : S_{G}(M) \to P$ be defined by $\phi(p)=px_0$. It is a continuous function, and $\phi|_{\cal M}$ is onto $P$.  Let us write $\phi'$ for the restriction of $\phi$ to ${\cal M}$. 

%Let $\phi = \phi_{x_{0}}: S_{G}(M) \to P$ be the unique surjective homomorphism taking $e$ to $x_{0}$ (see Fact \ref{universality of beta G} and the paragraph below it).  
%
%%%Krzys: I do not think that ux_0=x_0$. We can chose $x_0$ with this property, but it seems irrelevant in the proof. 
%Then $\phi(u) = x_{0}$ and 
%Then $\phi|_{\cal M}$ is onto $P$. 

%Let us write $\phi'$ for the restriction of $\phi$ to ${\cal M}$. For $p\in S_{G}(M)$ we write $px_{0}$ for $\phi_{x_{0}}(p)$.
%For $p,q\in S_{G}(M)$, we have $(pq)x_{0} = p(qx_{0}) (= \phi_{qx_{0}}(p))$. (All in Chapter 1 of \cite{Glasner}, where $pq$ means in the sense of the canonical semigroup structure on $S_{G}(M)$.)
 
Suppose for a contradiction that $P$ is non-trivial.\\
 
 \noindent
 {\bf Claim 2.} We can write ${\cal M}$ as the union of finitely many clopen subsets $V_{1},\dots,V_{n}$, such that for each $i$, $\phi'(V_{i})$ is a proper subset of $P$.\\[3mm] 
{\em Proof of Claim 2.} 
Easy using that $P$ is non-trivial, $\phi'$ is continuous and surjective, and ${\cal M}$ is a Stone space. Namely, first write $P$ as the union of finitely many proper open subsets $U_{1},\dots,U_{k}$. Let $V_{i}'$ be $\phi'^{-1}(U_{i})$. Then the $V_{i}'$ form a covering of $\cal M$. Each $V_{i}$ is open so a union of clopen subsets $V_{i,j}'$. By compactness, a finite union of some of the $V_{i,j}'$ cover $\cal M$ and this is the required covering. \hfill $\square$\\
 
%%%Krzys: I did not check the reference to Chernikov-Simon, so I am not sure about it. Is the number of the proposition in Ch-S 32 or rather 3.2? It is a pity that they have not circulated this paper yet, because the referee of our paper will not be able to check this (unless the referee will be Artem or Pierre). I added `continuous', because we use it later.
% Let now $\theta:{\cal M} \to G^*/{G^*}^{00}$ be the canonical, continuous, surjective homomorphism (it was defined on $u{\cal M}$ in Section \ref{proof of main theorem 1}, but the same definition works on ${\cal M}$). Now, Proposition 32 (in an early version) of \cite{ChSi} says that for every clopen subset $X$ of $\cal M$ (i.e. given by a formula over $M$), $\{g\in G^*/{G^*}^{00}: \theta^{-1}(g)$ meets both $X$ and $X^{c}\}$ has no interior. 
 %(This is strictly speaking over a saturated model but also works in our context.) 
Applying Claim 1 to all the $V_{i}$ from Claim 2, we obtain $g\in G^*/{G^*}^{00}$ and $i$ (without loss $i=0$) such that $\theta^{-1}(g)$ is contained in $V_{0}$.  Let $U_{0} = \phi'(V_{0})$. So $U_{0}$ is a proper subset of $P$. Let $p\in \theta^{-1}(g)$, and let $x = \phi'(p) = px_{0}$.  Now, let $y\in P\setminus U_{0}$. 
%%%Krzys: I have removed `(adapted to the current context)', because we use directly the version from Glasner. I also changed a little bit the next sentences.
By Proposition 3.2(2) of Chapter I of \cite{Gl}, 
%(adapted to the current context) 
there is some idempotent $v\in S_{G}(M)$ (even in a minimal ideal of $S_{G}(M)$) such that $vx=y$. As $v$ is an idempotent, $\hat{\theta}(v) =e{G^*}^{00}$ (where $\hat{\theta}$ is the obvious extension of $\theta$ to the function from $S_G(M)$ to $G^*/{G^*}^{00}$). Also, since ${\cal M}$ is a left ideal, $vp\in {\cal M}$. Therefore, $vp\in \theta^{-1}(g) \subseteq V_{0}$.  So $\phi'(vp) \in U_{0}$.  But $\phi'(vp) = (vp)x_{0} = v(px_{0}) = vx = y\notin U_{0}$, a contradiction. So the proof has been completed. \hfill $\square$\\

Let us discuss now the topology on the space of measures ${\cal M}(S_G(M))$.
A regular probability measure on $S_{G}(M)$ is the ``same thing'' as a finitely additive probability measure on the Boolean algebra $B_{G}(M)$ of clopens (i.e. a Keisler measure on $G$ over $M$). The collection of such finitely additive measures on $B_{G}(M)$ is a closed subset of $[0,1]^{B_{G}(M)}$ equipped with the product topology, yielding the structure of a compact Hausdorff space on ${\cal M}(S_{G}(M))$.  
On the other hand, ${\cal M}(S_{G}(M))$ can be identified (Riesz representation theorem) with a compact convex subset of the dual space to the space 
of continuous functions from $S_{G}(M)$ to $\R$ equipped with the weak $*$-topology. One can check that the two topologies coincide.

We will deduce below Proposition \ref{definability of the flow of measures}  from Theorem 2.7 of \cite{ChPiSi} which proves  definability of measures assuming definability of types, in the NIP environment. 
%We give some details of this, mainly recalling definitions and facts, as there could be some delicate points.

%%%Krzys: I changed (everywhere below) slightly the notation to be compatible with the rest of the paper. 

As in the case of groups (see Section \ref{Definable and externally definable context}), there are on the face of it two definitions of a definable map from a definable set to a compact space; the first, in a saturated environment, and the second in a not necessarily saturated environment. The definitions are set up to be compatible and we repeat the explanation.

Recall that $M$ denotes an arbitrary model, and $\C$ a monster model. For any set $X$ definable in $M$, by $X^*$ we denote its interpretation in $\C$.

%We identify notationally a definable set with its points in $\cal U$. 

\begin{definition}\label{definability for definable sets} 
Let $X$ by any set definable in $M$, and let $C$ be a compact Hausdorff space.\\
(i) A map $f:X^*\to C$ is said to be definable over $M$ if the premiage of any closed subset of $C$ is type-definable over $M$. (So this precisely means that $f$ is induced by a continuous map $\bar f: S_{X}(M) \to C$.)\\
(ii) A map $f:X\to C$ is said to be definable if for any two closed disjoint subsets $C_{1}, C_{2}$ of $C$, $f^{-1}[C_{1}]$ and $f^{-1}[C_{2}]$ are separated by a definable (with parameters from $M$ of course) subset of $X$.
\end{definition}

\begin{fact}\label{correspondence between definabilities} The map taking $f$ to $f|_X$ sets up a bijection between $M$-definable maps from $X^*$ to $C$ and definable maps from $X$ to $C$. 
\end{fact}

%We will consider Fact \ref{correspondence between definabilities} first in the context of definable measures and then in the context of definable actions.

As a prequel we consider definable types. For $\phi(x,y)$ an $L$-formula (or even $L_{M}$-formula), let $Y_{\phi}$ be the sort of the $y$ variable. Let $p(x)\in S_{x}(M)$. For $p$ to be definable means precisely that for each $\phi(x,y)$, the map $f_{p,\phi}$ from $Y_{\phi}(M)$ to $\{0,1\}$
such that $f_{p,\phi}(b) = 1$ iff $\phi(x,b)\in p$, is definable. A special case of Fact \ref{correspondence between definabilities} is that if $p\in S_{x}(M)$ is definable, then $p$ has a unique extension $p'\in S_{x}({\C})$ which is definable over $M$. 
%In fact,  the map taking $p$ to $p'$ establishes a homeomorphism between the set of definable types in $S_{x}(M)$ and the subset of $S_{x}(\C)$ consisting of global types definable over $M$. 

Now for measures. Let $X$ be an $M$-definable set, and $\mu$ a Keisler measure on $X$ over $M$, namely $\mu\in{\cal M}(S_{X}(M))$ according to earlier identifications.  To say that $\mu$ is definable needs to be understood in terms of Definition \ref{definability for definable sets}(ii), as we are in a not necessarily saturated environment. So $\mu$ is {\em definable} if for all $\phi(x,y)$, the map $f_{\mu,\phi}:Y_{\phi}(M)\to [0,1]$ taking $b$ to $\mu(\phi(x,b))$  is definable in the sense of Definition \ref{definability for definable sets}(ii). 
%
%%%Krzys: The fact that the unique definable over $M$ extension is really a measure requires an easy proof, but I guess we can skip it. I replaced `means' by `implies', and then `taking' by `from'.
%
And Fact \ref{correspondence between definabilities} implies that such a definable $\mu$ has a unique extension to a global Keisler measure $\mu'$ on $X$ which is definable over $M$ in the sense of Definition \ref{definability for definable sets}(i): namely for each $\phi$ the map from $Y_{\phi}^*$ to $[0,1]$ taking $b$ to $\mu'(\phi(x,b))$ is definable over $M$ in the sense of Definition \ref{definability for definable sets}(i). 
%The mapping taking $\mu$ to $\mu'$ again establishes a homeomorphism between the set of definable Keisler measures on $X$ over $M$ and the set of global Keisler measures on $X$ %which are definable over $M$ (where the topologies are as described above). 

Now for definable actions. We are in the context of a group $G$ definable in $M$. As was recalled in Definition \ref{definition of definable flows}, an action of $G$ on a compact space $C$ by homeomorphisms is said to be definable if for each $x\in C$ the map $f_{x}:G \to C$ defined by $f_{x}(g) = gx$ is definable in the sense of Definition \ref{definability for definable sets}(ii).
We will be interested in proving that the action of $G$ on the compact space ${\cal M}(S_{G}(M))$ is definable.\\ 

 %By Fact \ref{correspondence between definabilities}, this precisely means that for each $x\in C$ there is map from $G^*$ to $C$ which is definable over $M$ in the sense of Definition %\ref{definability for definable sets}(i) and whose restriction to $G$ is $f_{x}$. (But one should be careful: it does {\em not}  mean that the action of $G$ on $C$ is the restriction to $G$ of a %suitable action of $G^*$ on $C$.)

\noindent
{\em Proof of Proposition \ref{definability of the flow of measures}.} 
Theorem 2.7 from \cite{ChPiSi} shows (using the VC-theorem) that if $\mu$ is a Keisler measure on $G$ over $M$, then $\mu$ is definable in the sense of Definition \ref{definability for definable sets}(ii).  We want to deduce from this that the action of $G$ on ${\cal M}(S_{G}(M))$ is definable. Fix $\mu\in {\cal M}(S_{G}(M))$, and we have to prove that the map, which we will call $f_{\mu}$,  taking $g\in G$ to $g\mu$ is definable.  So we fix disjoint closed subsets $D_{1}$ and $D_{2}$ of ${\cal M}(S_{G}(M))$ and we have to show that $f_{\mu}^{-1}[D_{1}]$ and $f_{\mu}^{-1}[D_{2}]$ are separated by a definable subset of $G$.  By the description above of the topology on ${\cal M}(S_{G}(M))$, we may assume that $D_{1}$ is of the form $\{\rho\in {\cal M}(S_{G}(M)): \rho(X_{1})\in C_{1}\; \mbox{and}\; \rho(X_{2})\in C_{2}\; \mbox{and}\; \dots \; \mbox{and} \;\rho(X_{k})\in C_{k}\}$, where $X_{1}, \dots,X_{k}$ are definable (over $M$) subsets of $G$ and $C_{1},\dots, C_{k}$ are closed subsets of $[0,1]$.  And similarly for $D_{2}$. 

Now let $\mu'$ be the unique extension of $\mu$ to a global Keisler measure which is definable over $M$ (given by  Fact \ref{correspondence between definabilities} as discussed above). 
For $g\in G^{*}$, let $g\mu'|_M$ denote the restriction of the global measure $g\mu'$ to $M$. So $g\mu'|_M \in  {\cal M}(S_{G}(M))$. 
Let $Y_{1} = \{g\in G^{*}: g\mu'|_M \in D_{1}\}$, which is precisely $\{g\in G^{*}: \mu'(g^{-1}X_{1})\in C_{1}\; \mbox{and}\; \dots \; \mbox{and}\; \mu'(g^{-1}X_{k})\in C_{k}\}$, which is type-definable over $M$, by definability of $\mu'$.  Similarly $Y_{2} = \{g\in G^{*}:g\mu'|_M \in D_{2}\}$ is type-definable over $M$.  As $Y_{1}$ and $Y_{2}$ are disjoint, we obtain by compactness an $M$-definable subset $Y$ of $G$ such that $Y^*$ separates them. Then clearly $Y$ separates $f_{\mu}^{-1}[D_{1}]$ and $f^{-1}_{\mu}[D_{2}]$ as required. \hfill $\square$\\

\noindent
{\em Proof of Proposition \ref{strong amenability implies amenability}.} 
%%%Krzys: I added the next sentence.
We follow the lines of the proof that strong amenability implies amenability (see \cite[Chapter III, Theorem 3.1]{Gl}). Recall that definable strong amenability of $G$ means that there is no nontrivial definable minimal proximal $G$-flow. 
We are assuming definable strong amenability of $G$ and want to deduce definable amenability, namely the existence of a $G$-invariant Keisler measure on $G$ over $M$.  
%%%Krzys: I removed the reference to Proposition 6.3.
By assumption, the natural action of $G$ on ${\cal M}(S_{G}(M))$ by homeomorphisms is definable. 
By Theorem 2.3 of Chapter III of \cite{Gl}, there is a minimal $G$-subflow $X$ of ${\cal M}(S_{G}(M))$ which is strongly proximal. Strong proximality means that the action of $G$ on ${\cal M}(X)$ is proximal. 
%%%Krzys: I removed `closed invariant', as subflows are such by definition.
Now, $X$ is a definable $G$-flow (as a subflow of the definable $G$-flow ${\cal M}(S_{G}(M))$). 
%%%Krzys: I did not check that $M(X)$ is also definable.
One can show that the action of $G$ on ${\cal M}(X)$ is also definable, but we will not need it. Now, $X$ embeds as a $G$-flow, homeomorphically into ${\cal M}(X)$: map $x\in X$ to the measure concentrating on $x$. Hence, $X$ is a proximal $G$-flow. As $X$ is minimal and $G$ is definably strongly amenable, $X$ is a singleton $\{x\}$, and so $x$ is a $G$-invariant measure on $S_{G}(M)$, yielding definable amenability of $G$. \hfill $\square$

\section{A characterization of the externally definable Bohr compactification}

Our goal here is to give a description (in the spirit of Fact \ref{Bohr compactification as quotient}) of the externally definable Bohr compactification of a given group as a certain quotient by a connected component. However, in order to do that, we will have to work in a certain subgroup (which is invariant but not necessarily type-definable) of the given group. We will also get appropriate counterparts of Theorem \ref{main theorem 1} and Corollary \ref{main corollary 3}.

For the rest of this section, let $H=\langle S_{G,M}(N)(\C)\rangle \leq G^*$. So $H$ is an $N$-invariant subgroup. (This group depends on the choice of $N$, so it should be denoted by something like $H_N$, but we denote it by $H$ for simplicity.) Let $H^{000}_N$ be the smallest $N$-invariant subgroup of $H$ of bounded index. As usual, one easily shows that $H^{000}_N$ exists and is a normal subgroup of $H$, and one gets an absolute bound on the index $[H:H^{000}_N]$ (using Erd\H{o}s-Rado theorem). The proof of the next proposition is written in such way that it shows this, too.

\begin{proposition}\label{H^000_N}
$H^{000}_N = \langle a^{-1}b: \tp(a/N)=\tp(b/N) \in S_{G,M}(N) \rangle$.
\end{proposition}
{\em Proof.} Let $H_0$ be the right hand side of this equality.\\[1mm]
%The inclusion $(\subseteq)$ is clear since $H_0$ is $N$-invariant and $[H:H_0] \leq  |S_{G,M}(N)| \cdot \aleph_0$. 
$(\subseteq)$ First, we show that $H_0$ is normal in $H$. For this, it is enough to check that for any $a,a',b$ such that $\tp(a/N)=\tp(a'/N) \in S_{G,M}(N)$ and $\tp(b/N) \in S_{G,M}(N)$ one has that $b^{-1}(a^{-1}a')b \in H_0$. We can find $b'$ realizing the unique extension of $\tp(b/N)$ to a complete type over $N,a,a'$ finitely satisfiable in $M$. Then $b'^{-1}b \in H_0$, so it is enough to show that $b'^{-1}(a^{-1}a')b' \in H_0$. By the fact that $\tp(b'/N,a,a')$ is finitely satisfiable in $M$ and $\tp(a/N)=\tp(a'/N)$, we easily get that $\tp(a/N,b')=\tp(a'/N,b')$. Hence, $\tp(b'^{-1}ab'/N)=\tp(b'^{-1}a'b'/N)$. Using the fact that $\tp(b'/N,a)$ and $\tp(a/N)$ are finitely satisfiable in $M$, we easily get that $\tp(b'^{-1}ab'/N) \in S_{G,M}(N)$. Thus, 
$b'^{-1}(a^{-1}a')b'=(b'^{-1}ab')^{-1}(b'^{-1}a'b') \in H_0$.

By the normality of $H_0$, one easily gets that $[H:H_0] \leq  |S_{G,M}(N)| \cdot \aleph_0$. Since clearly $H_0$ is $N$-invariant, we conclude that $H_0$ is an $N$-invariant normal subgroup of $H$ of bounded index.\\[1mm]
$(\supseteq)$ Let $L$ be any $N$-invariant subgroup of $H$ of bounded index. Consider any $a,b$ such that $\tp(a/N)=\tp(b/N) \in S_{G,M}(N)$. Then there is $c \models \tp(a/N)$ such that both pairs $(a,c)$ and $(b,c)$ can be extended to arbitrarily long infinite sequences indiscernible over $N$. We will show that $a^{-1}c \in L$ and $c^{-1}b\in L$, which implies $a^{-1}b \in L$, finishing the proof.

Let $(a_i)_i$ be an $N$-indiscernible sequence extending $(a,c)$ of length greater that $[H:L]$. Note that all $a_i$'s are in $H$. If $a^{-1}c\notin L$, then $a_i^{-1}a_j \notin L$ for all $i<j$, which contradicts the choice of the length of the sequence  $(a_i)_i$. \hfill $\square$

\begin{remark}\label{map pi}
The map $\pi: S_{G,M}(N)(\C) \to H/H^{000}_N$ given by $\pi(a)=aH^{000}_N$ is onto.
\end{remark}
{\em Proof.} 
It is enough to show that whenever $\tp(a/N),\tp(b/N) \in S_{G,M}(N)$, one has $abH^{000}_N \cap S_{G,M}(N)(\C) \ne \emptyset$. Take $a'$ such that $\tp(a'/N,b)$ is the unique extension of $\tp(a/N)$ to a complete type over $N,b$ which is finitely satisfiable in $M$. Then $\tp(a'b/N) \in S_{G,M}(N)$, and, by Proposition \ref{H^000_N}, $(a'b)^{-1}ab=b^{-1}a'^{-1}ab \in b^{-1}H^{000}_Nb=H^{000}_N$, so $a'b \in abH^{000}_N$. \hfill $\square$\\

Now, we introduce the logic topology on $H/H^{000}_N$.

\begin{definition}\label{logic topology}
We say that $D\subseteq H/H^{000}_N$ is closed if $\pi^{-1}[D]$ is a type-definable subset of $S_{G,M}(N)(\C)$ (equivalently, if it is type-definable over $N$), where $\pi$ is the map defined in Remark \ref{map pi}.
\end{definition}

Let $H^{00}_N\leq H$  be the preimage of the closure of the neutral element of $H/H^{000}_N$ by the quotient map from $H$ to $H/H^{000}_N$.

\begin{proposition}
$H^{00}_N$ is a normal subgroup of $H$.
\end{proposition}
{\em Proof.}
First we show that for any $a \in H$ the function $f_a: H/H^{000}_N \to H/H^{000}_N$ given by $f_a(xH^{000}_N)=axH^{000}_N$ is continuous (i.e. the group operation in $H/H^{000}_N$ is continuous on the right). By Remark \ref{map pi}, we can choose $a' \in aH^{000}_N$ such that $\tp(a'/N) \in S_{G,M}(N)$. Consider any closed set $D \subseteq H/H^{000}_N$. Then 
$$\pi^{-1}[f_a^{-1}[D]]= \{ b \in S_{G,M}(N)(\C) : (\exists b' \equiv_N b)( a'b' \in \pi^{-1}[D])\}$$ 
is type-definable, so $f_a^{-1}[D]$ is closed.

Since $f_{a^{-1}}$ is the inverse of $f_a$, we get that $f_a$ is a homeomorphism. Similarly, one gets that inversion and conjugation by any element are homeomorphisms of $H/H^{000}_N$.

So, in order to finish the proof, it is enough to show that $H^{00}_N$ is closed under multiplication. Since for any $a\in H$, $f_a$ is a homeomorphism, we get that 
$\cl(aH^{000}_N)=a\cl(eH^{000}_N)$. Now, consider any $a,b \in H^{00}_N$. We have $abH^{000}_N \in a\cl(eH^{000}_N) =\cl(aH^{000}_N) \subseteq \cl(eH^{000}_N) = H^{00}_N/H^{000}_N$, so $ab \in H^{00}_N$. \hfill $\square$\\

Let $\bar \pi: S_{G,M}(N)(\C) \to H/H^{00}_N$  be the map given by $\bar \pi(a)=aH^{00}_N$. By Remark \ref{map pi}, $\bar \pi$ is onto.
We introduce the  logic topology on $H/H^{00}_N$ as in Definition \ref{logic topology}.

%By routine but not completely obvious arguments  (using the facts that the map $\pi$ factors through the natural map $r: S_{G,M}(N)(\C) \to S_{G,M}(N)$ and that $S_{G,M}(N)$ is closed in $S_G(N)$), one can show the following proposition.

\begin{proposition}\label{H/H^00 is compact}
i) $H/H^{00}_N$ is a compact Hausdorff topological group.\\
ii) $H/H^{000}_N$ is a compact topological group (which is not necessarily Hausdorff).
\end{proposition}
{\em Proof.} 
(i) Compactness of $H/H^{00}_N$ is clear from the definition of the logic topology and the fact that $\bar \pi$ is onto. 
%and that $S_{G,M}(N)$ is closed.

Now we will show that $H/H^{00}_N$ is Hausdorff. Consider any $aH^{00}_N \in H/H^{00}_N$ and choose  $a' \in aH^{00}_N \cap S_{G,M}(N)(\C)$. One easily gets
$$\bar \pi^{-1}(aH^{00}_N)=\left\{ b \in S_{G,M}(N)(\C) : (\exists b' \equiv_N b) (a'^{-1}b' \in \pi^{-1}[H^{00}_N/H^{000}_N]) \right\}.$$
$\pi^{-1}[H^{00}_N/H^{000}_N]$ is an $N$-type-definable set, so it is defined by a collection of formulas $\theta_i(x)$, $i \in I$, over $N$ closed under conjunction. Choose also a collection of formulas $\varphi_j(x,y)$, $j \in J$, over $N$ which defines the condition $x \equiv_N y$ and which is closed under conjunction. For any $i \in I$ and $j \in J$ consider
$$O_{a',i,j}:= \left\{cH^{00}_N: \bar \pi^{-1}(cH^{00}_N) \subseteq \{ b \in S_{G,M}(N)(\C) : (\exists b') (\varphi_j(b',b) \wedge \theta_i (a'^{-1}b'))\}\right\}.$$
Clearly $aH^{00}_N \in O_{a',i,j}$, and one can check that $O_{a',i,j}$ is open in the logic topology on $H/H^{00}_N$.

Consider any $a_1H^{00}_N \ne a_2H^{00}_N$ in $H/H^{00}_N$. Take $a_1' \in a_1H^{00}_N \cap S_{G,M}(N)(\C)$ and $a_2' \in a_2H^{00}_N \cap S_{G,M}(N)(\C)$. We claim that there are $i \in I$ and $j \in J$ such that $O_{a_1',i,j} \cap O_{a_2',i,j} = \emptyset$, which of course witnesses that $H/H^{00}_N$ is Hausdorff. If it is not the case, then, by compactness, there are $b',b'' \in S_{G,M}(N)(\C)$ such that $b' \equiv_N b''$ and $a_1'^{-1} b', a_2'^{-1}b''  \in  \pi^{-1}[H^{00}_N/H^{000}_N]$, so $a_1'^{-1} b'b''^{-1} a_2' \in H^{00}_N$, so $a_1'H^{00}_N=a_2'H^{00}_N$, a contradiction.

The fact that group inversion in $H/H^{00}_N$ is continuous is very easy. Let us check that multiplication is continuous. The map $\bar \pi$ factors through the natural map $r: S_{G,M}(N)(\C) \to S_{G,M}(N)$, i.e., there is a continuous map $p: S_{G,M}(N) \to H/H^{00}_N$ such that $\bar \pi=p\circ r$.  Then $\bar \pi \times \bar \pi =(p \times p) \circ (r \times r)$. Denote the group operation in $H/H^{00}_N$ by $\cdot$. 

Take any closed subset $D \subseteq H/H^{00}_N$, and put $S:=(\cdot \circ (\bar \pi \times \bar \pi))^{-1}[D]$. Then
$$S=\left\{ (a,b) \in S_{G,M}(N)(\C) \times  S_{G,M}(N)(\C) : (\exists a' \equiv_N a)(\exists b'\equiv_N b) (a'b' \in \bar \pi^{-1}[D])\right\},$$
which is of course type-definable. Thus, $(r\times r) [S]$ is closed, and so $((p \times p) \circ (r \times r)) [S]$ is closed, too (as $S_{G,M}(N)$ is compact, $H/H^{00}_N$ is Hausdorff and $p \times p$ is continuous). But the last set is exactly $\cdot^{-1}[D]$, so we are done.\\[1mm]
(ii) is left as an exercise. \hfill $\square$

\begin{remark}\label{pi[G] is dense}
i) $\pi[G]$ is a dense subgroup of $H/H^{000}_N$.\\
ii) $\bar\pi [G]$ is a dense subgroup of $H/H^{00}_N$.
\end{remark}
{\em Proof.} 
It is enough to show (i). Take the induced map $h: S_{G,M}(N) \to H/H^{000}_N$. Consider any open subset $U$ of $H/H^{000}_N$. Then, there is a formula $\varphi(x)$ with parameters from $N$ such that $\varphi(G) \ne \emptyset$ and $[\varphi(x)] \cap S_{G,M}(N) \subseteq h^{-1}[U]$. Take $g \in \varphi(G)$. Then $gH^{000}_N \in U$. \hfill $\square$\\

By Proposition \ref{H/H^00 is compact}(i), Remark \ref{pi[G] is dense}(ii), the definition of the logic topology on $H/H^{00}_N$ and by point 2)(ii) of Lemma \ref{prolongation}, we get

\begin{corollary}
$H/H^{00}_N$ is an externally definable compactification of $G$.
\end{corollary}

Now, we can prove the promised description of the externally definable Bohr compactification of $G$.

\begin{proposition}
$H/H^{00}_N$ is the externally definable Bohr compactification of $G$.
\end{proposition}
{\em Proof.}
Let $C$ be an externally definable compactification of $G$; so there is an externally definable homomorphism $f: G \to C$, and Lemma \ref{prolongation} yields the externally definable map $f^*:S_{G,M}(N)(\C) \to C$ given by  $f^*(a)=\bigcap_{\varphi \in \tp(a/N)} \cl(f[\varphi(M)])$.\\[3mm]
{\bf Claim} If $a,b \in S_{G,M}(N)(\C)$ and $aH^{000}_N=bH^{000}_N$, then $f^*(a)=f^*(b)$.\\[3mm]
{\em Proof of the claim.} By Remark \ref{H^000_N}, there are some $a_1,\dots,a_n, a_1',\dots,a_n' \in H$ satisfying $\tp(a_i/N)=\tp(a_i'/N) \in S_{G,M}(N)$ and such that $a=b(a_1^{-1}a_1')\cdot \ldots \cdot (a_n^{-1}a_n')$. Since $a,b$ and all $a_i$'s and $a_i'$'s are in the domain of $f^*$, Lemma \ref{partial homomorphism} gives us $f^*(a)=f^*(b) \prod_i f^*(a_i^{-1})f^*(a_i')=f^*(b) \prod_i f^*(a_i)^{-1}f^*(a_i)=f^*(b)$, which completes the proof of the claim. \hfill $\square$\\

By the claim and Remark \ref{map pi}, we get a well-defined map $\Phi: H/H^{000}_N \to C$ given by 
$$\Phi(aH^{000}_N) =f^{*}(a')$$ 
for any $a' \in aH^{000}_N \cap S_{G,M}(N)(\C)$.

Let us check that $\Phi$ is a homomorphism. Consider any $a,b \in H$. Take $a',b'\in S_{G,M}(N)(\C)$ in the cosets $aH^{000}_N, bH^{000}_N$, respectively, so that $\tp(a'/N,b')$ is finitely satisfiable in $M$. Then $a'b' \in abH^{000}_N \cap  S_{G,M}(N)(\C)$. By Lemma \ref{partial homomorphism}, we get
$\Phi(aH^{000}_N \cdot bH^{000}_N)=\Phi(a'b'H^{000}_N)=f^*(a'b')=f^*(a')f^*(b')=\Phi(aH^{000}_N)\Phi(bH^{000}_N)$.

The fact that $\Phi$ is continuous is clear, as for a closed subset $D$ of $C$ one has $\pi^{-1}[\Phi^{-1}[D]]={f^*}^{-1}[D]$ is a type-definable subset of $S_{G,M}(N)(\C)$.

Thus, $\ker(\Phi)$ is a closed subgroup of $H/H^{000}_N$, so it contains $H^{00}_N/H^{000}_N$. This yields a natural continuous homomorphism from $H/H^{00}_N$ to $(H/H^{000}_N)/\ker(\Phi) $ commuting with the natural homomorphisms from $G$. Since $H/H^{000}_N$ is compact and $C$ is compact Hausdorff, we have that $(H/H^{000}_N)/\ker(\Phi) $ is topologically isomorphic to $C$ via the natural isomorphism also commuting with the homomorphisms from $G$ (note that for $g \in G$, one has $\Phi(gH^{000}_N)=f^*(g)=f(g)$). Thus, we have found a continuous homomorphism from $H/H^{00}_N$ to $C$ commuting with the homomorphisms from $G$. \hfill $\square$\\

The proof of Theorem \ref{main theorem 1} goes through, and we get the following result.

\begin{proposition}
Let ${\cal M}$ be a minimal ideal in $S_{G,ext}(M)=S_{G,M}(N)$ and let $u\in {\cal M}$ be an idempotent. Then, the function
$f: u{\cal M} \to  H/H^{000}_N$ given by $\tp(a/N) \mapsto aH^{000}_N$ is a well-defined epimorphism. Equip $u{\cal M}$ with the $\tau$-topology. Then
\begin{enumerate}
\item $f$ is continuous.
\item $H(u{\cal M}) \leq f^{-1}[H^{000}_N]$.
\item The formula $\tp(a/N)H(u{\cal M}) \mapsto aH^{000}_N$ yields a well-defined continuous epimorphism from $u{\cal M}/H(u{\cal M})$ to $H/H^{000}_N$.
\end{enumerate}
\end{proposition}

Thus, as in the case of Corollary \ref{main corollary 3}, we get

\begin{corollary}\label{cor 2''}
Suppose $G$ is externally definably strongly amenable (i.e., the only minimal externally definable proximal flow is the trivial one). Then the natural epimorphism $\zeta: u{\cal M}/H(u{\cal M}) \to H/H^{00}_N$ is an isomorphism, so $H^{000}_N=H^{00}_N$.
\end{corollary}

One can also check that the proof of Theorem \ref{main theorem 2} goes through in our current context, yielding the obvious counterpart of this theorem.

All of this leads to various natural questions, some of which could probably be answered negatively by giving suitable counter-examples. We discuss only some of them. 

Since $H/H^{00}_N$ is the externally definable Bohr compactification of $G$, it does not depend on the choice of $N$. Is $H/H^{000}_N$ also independent of the choice of $N$ (recall that $H$ depends on $N$)? We think that this should be true. If $N_2 \succ N_1\succ M$ are $|M|^+$-saturated, then there is a continuous epimorphism from the quotient obtained for $N_2$ to the one obtained for $N_1$, but it is not clear whether it is injective.

There is also a question whether the externally definable Bohr compactification coincides with the definable Bohr compactification. Probably not. It would be also interesting to describe the externally definable Bohr compactification in terms of connected components of $G^*$ (not of $H$).  For example, $G^*/{G^*}^{00}_M$ as well as $\cl (G/{G^*}^{00}_N)$ (where $N\succ M$ is $|M|^+$-saturated) are externally definable compactifications.
Is any of them the externally definable Bohr compactification? Is the natural continuous epimorphism from $\cl (G/{G^*}^{00}_N)$ to $G^*/{G^*}^{00}_M$ an isomorphism? Is  $\cl (G/{G^*}^{00}_N)$ independent of the choice of $N$? 
There is a natural continuous epimorphism from $H/H^{00}_N$ to $\cl (G/{G^*}^{00}_N)$ and also to $G^*/{G^*}^{00}_M$. Is any of them an isomorphism? They are isomorphisms when all types in $S_G(M)$ are definable, as then the externally definable Bohr compactification coincides with the definable Bohr compactification.
%I have a characterizarion when it is true from which one can see directly that the second epimorphism (and so the first one, too) is injective under the assumption that all types in $S_G(M)$ are definable (but this remark follows immediately from the fact that under this assumption, $G/{G^*}^{00}_M$ is the externally definable Bohr compactification of $G$).

\vspace{2mm}

\noindent
{\bf Addresses:}\\[2mm]
Krzysztof Krupi\'nski\\[1mm]
Instytut Matematyczny, Uniwersytet Wroc\l awski\\
pl. Grunwaldzki 2/4\\ 
50-384 Wroc\l aw, Poland\\[1mm]
{\bf e-mail}: kkrup@math.uni.wroc.pl\\[3mm]
Anand Pillay\\[1mm]
Department of Mathematics, University of Notre Dame\\
281 Hurley Hall\\
Notre Dame, IN 46556, USA\\[1mm] 
{\bf e-mail}: apillay@nd.edu

\end{document}